\documentclass[a4paper, 11pt, twoside]{article}
\usepackage{amsmath}
\usepackage{amscd}
\usepackage{amssymb}
\usepackage{amsthm}
\usepackage{graphics}
\usepackage{indentfirst}
\usepackage{latexsym}
\usepackage{amsfonts}
\usepackage{multicol}
\usepackage{stmaryrd}
\usepackage{array}
\usepackage{youngtab}
\usepackage{pxfonts}
\usepackage{verbatim}
\usepackage{pb-diagram}

\newtheoremstyle{mattthm}{}{}{\itshape}{}{\bfseries}{.}{ }{}

\theoremstyle{mattthm}
\newtheorem{lemma}{Lemma}[section]
\newtheorem{propn}[lemma]{Proposition}
\newtheorem{thm}[lemma]{Theorem}
\newtheorem{cory}[lemma]{Corollary}

\newtheoremstyle{mattdef}{}{}{}{}{\bfseries}{.}{ }{}

\theoremstyle{mattdef}

\newtheorem*{defn}{Definition}
\newtheorem*{eg}{Example}
\newtheorem*{egs}{Examples}
\newtheorem*{rmks}{Remarks}
\newtheorem*{acks}{Acknowledgements}

\begin{document}

\newenvironment{pf}{\noindent\textbf{Proof.}}{\hfill\qedsymbol\newline}
\newenvironment{pfof}[1]{\vspace{\topsep}\noindent\textbf{Proof of {#1}.}}{\hfill\qedsymbol\newline}
\newenvironment{pfenum}{\noindent\textbf{Proof.}\indent\begin{enumerate}\vspace{-\topsep}}{\hfill \qedsymbol\end{enumerate}\vspace{-\topsep}}
\newenvironment{pfnb}{\noindent\textbf{Proof.}}{\newline}
\newenvironment{pfofnb}[1]{\vspace{\topsep}\noindent\textbf{Proof of {#1}.}}{\newline}

\hyphenation{mod-ules}
\newcommand\co{\geqslant}
\newcommand\cso{>}
\newcommand\nco{\ngeqslant}
\newcommand\ar[1]{\operatorname{AR}_{#1}}
\newcommand\La\Lambda
\newcommand\bi{\ol\imath}
\newcommand\bj{\ol\jmath}
\newcommand\wt{\operatorname{wt}}
\newcommand\lra\longrightarrow
\newcommand\calc{\mathcal{C}}
\newcommand\g\gtrdot
\newcommand\gc{\operatorname{GC}}
\newcommand\cc{\operatorname{C}}
\newcommand\A{\operatorname{A}}
\newcommand\R{\operatorname{R}}
\newcommand\bsm{\begin{smallmatrix}}
\newcommand\esm{\end{smallmatrix}}
\newcommand{\rt}[1]{\rotatebox{90}{$#1$}}
\newcommand\la\lambda
\newcommand\ep\epsilon
\newcommand{\ol}{\overline}
\newcommand{\lan}{\langle}
\newcommand{\ran}{\rangle}
\newcommand\partn{\mathcal{P}}
\newcommand{\thmlc}[3]{\textup{\textbf{(\!\! #1 \cite[#3]{#2})}}}
\newcommand{\sss}{\mathfrak{S}_}
\newcommand{\dom}{\trianglerighteqslant}
\newcommand{\doms}{\vartriangleright}
\newcommand{\ndom}{\ntrianglerighteqslant}
\newcommand{\ndoms}{\not\vartriangleright}
\newcommand{\domby}{\trianglelefteqslant}
\newcommand{\domsby}{\vartriangleleft}
\newcommand{\ndomby}{\ntrianglelefteqslant}
\newcommand{\ndomsby}{\not\vartriangleleft}
\newcommand{\subs}[1]{\subsection{#1}}
\newcommand{\subsubs}{\subsubsection}
\newcommand{\nin}{\notin}
\newcommand{\nchar}{\operatorname{char}}
\newcommand{\thmcite}[2]{\textup{\textbf{\cite[#2]{#1}}}\ }
\newcommand\znz{\mathbb{Z}/n\mathbb{Z}}
\newcommand\zepz{\mathbb{Z}/(e+1)\mathbb{Z}}
\newcommand{\bbf}{\mathbb{F}}
\newcommand{\bbc}{\mathbb{C}}
\newcommand{\bbn}{\mathbb{N}}
\newcommand{\bbq}{\mathbb{Q}}
\newcommand{\bbz}{\mathbb{Z}}
\newcommand\zo{\bbn_0}
\newcommand{\gs}{\geqslant}
\newcommand{\ls}{\leqslant}
\newcommand\dw{^\triangle}
\newcommand\wod{^\triangledown}
\newcommand{\hhh}{\mathcal{H}_}
\newcommand{\sect}[1]{\section{#1}}
\newcommand{\fkf}{\mathfrak{f}}
\newcommand{\fff}{\mathfrak{F}}
\newcommand\cf{\mathcal{F}}
\newcommand\fkg{\mathfrak{g}}
\newcommand\fkh{\mathfrak{h}}
\newcommand\fkk{\mathfrak{k}}
\newcommand\fkl{\mathfrak{l}}
\newcommand\fkm{\mathfrak{m}}
\newcommand\fkn{\mathfrak{n}}
\newcommand\sx{x}
\newcommand\bra[1]{|#1\ran}
\newcommand\arb[1]{\widehat{\bra{#1}}}
\newcommand\foc[1]{\mathcal{F}_{#1}}
\newcommand{\clam}{\begin{description}\item[\hspace{\leftmargin}Claim.]}
\newcommand{\prof}{\item[\hspace{\leftmargin}Proof.]}
\newcommand{\malc}{\end{description}}
\newcommand\ppmod[1]{\ (\operatorname{mod}\ #1)}
\newcommand\wed\wedge
\newcommand\wede\barwedge
\newcommand\uu[1]{\,\begin{array}{|@{\,}c@{\,}|}\hline #1\\\hline\end{array}\,}
\newcommand{\ux}{\operatorname{sz}}
\newcommand\erim{\operatorname{rim}}
\newcommand\mire{\operatorname{rim}'}
\newcommand\mmod{\ \operatorname{Mod}}
\newcommand\cgs\succcurlyeq
\newcommand\cls\preccurlyeq
\newcommand\cg\succ
\newcommand\cl\prec
\newcommand\ee{\tilde{e}_}
\newcommand\ff{\tilde{f}_}
\newcommand\vh{\Gamma}
\newcommand\fsl{\mathfrak{sl}}
\newcommand\hsl{\widehat{\mathfrak{sl}}}
\newcommand\reg[1]{\mathcal{R}_{#1}}
\newcommand\regg[1]{\mathcal{R}_{#1}}
\newcommand\zez{\bbz/e\bbz}
\newcommand\cla{\mathcal{C}}
\newcommand\cld{\mathcal{D}}
\newcommand\claup{\mathcal{C}^{\wedge}}
\newcommand\cladw{\mathcal{C}^{\vee}}
\newcommand\rest[1]{\operatorname{Rest}_{#1}}
\newcommand{\ja}{\put(10,-5){\vector(2,-1){40}}\put(30,-13.5){\small$j$}}
\newcommand{\ia}{\put(-10,-5){\vector(-2,-1){40}}\put(-30,-13.5){\small$i$}}
\newcommand{\jaa}{\put(7.5,-5){\vector(3,-2){45}}\put(30,-18){\small$j$}}
\newcommand{\iaa}{\put(-7.5,-5){\vector(-3,-2){45}}\put(-30,-18){\small$i$}}
\newcommand{\vt}[1]{\put(0,-2.2){\makebox[0pt]{$(#1)$}}}
\newcommand\ny[2]{\node[#1]{\yng(#2)}}
\newcommand\are{\arrow{e,t}}
\newcommand\ane{\arrow{ne,t}}
\newcommand\ase{\arrow{se,t}}
\newlength\spacer
\setlength\spacer{12pt}

\title{Partition models for the crystal of the basic $U_q(\widehat{\mathfrak{sl}}_n)$-module}
\author{Matthew Fayers\\\normalsize Queen Mary, University of London, Mile End Road, London E1 4NS, U.K.\\\texttt{\normalsize m.fayers@qmul.ac.uk}}
\date{}
\maketitle
\begin{center}
\end{center}
\markboth{Matthew Fayers}{Partition models for the crystal of the basic $U_q(\widehat{\mathfrak{sl}}_n)$-module}
\pagestyle{myheadings}

\begin{abstract}
For each $n\gs3$, we construct an uncountable family of models of the crystal of the basic $U_q(\hsl_n)$-module.  These models are all based on partitions, and include the usual $n$-regular and $n$-restricted models, as well as Berg's ladder crystal, as special cases.
\end{abstract}

\sect{Introduction}

Let $n\gs2$, and let $\hsl_n$ denote the Kac--Moody algebra of type $A^{(1)}_{n-1}$, and $U_q(\hsl_n)$ its quantised enveloping algebra.  A key object in the combinatorial representation theory of a quantised Kac--Moody algebra $U$ is a \emph{crystal basis} for a module $M$, and the associated \emph{crystal graph}.  When $U=U_q(\hsl_n)$ and $M$ is the irreducible highest-weight module with highest weight $\La_0$, this crystal graph plays a r\^ole in the representation theory of Iwahori--Hecke algebras of type $A$ at an $n$th root of unity: there is a correspondence between the vertices of this crystal graph and simple modules for these Hecke algebras, with the arrows in the crystal graph corresponding to certain functors refining induction and restriction between these algebras.

There is particular interest in modelling crystal graphs, i.e.\ realising their vertex sets as sets of well-understood combinatorial objects.  Of particular prominence in the case of $\hsl_n$-crystals is the modelling of crystals of highest-weight representations by multipartitions (where the number of components in each multipartition equals the level of the highest weight).  In this paper we restrict attention to the better-understood case of level $1$, where crystals are modelled by partitions.  There are two well-established models of the highest-weight crystal with highest weight $\La_0$, in terms of $n$-regular and $n$-restricted partitions.  These models come directly from \emph{Fock space} models for (modules containing) the corresponding highest-weight module $V(\La_0)$.

More recently, Berg \cite{berg} found a realisation of $B(\La_0)$ in terms of a new class of partitions, in the case where $n\gs3$.  Berg's construction is purely combinatorial -- his proof involves an explicit isomorphism of crystals between his crystal and the $n$-regular crystal -- and it is unclear whether there is more algebraic structure behind the scenes.

In this paper we extend Berg's work, by showing that his crystal model is one of an uncountable family of models of the crystal $B(\La_0)$, for each $n\gs3$.  Again, our construction and proof are purely combinatorial.  The proof is rather awkward -- it is not clear how to give an explicit isomorphism from one of our crystals to a known model of $B(\La_0)$, and so we use a result of Stembridge, which enables recognition of crystals of highest-weight modules for simply-laced algebras in terms of local properties.  In fact, it is too difficult for this author even to use Stembridge's results directly, and we use an alternative route due to Danilov, Karzanov and Koshevoy.  In some ways this proof is unsatisfactory; it would be nice to have an algebraic understanding -- perhaps via an analogue of the Fock space -- for why these models exist.

We now indicate the layout of this paper.  In the next section we give enough definitions to enable a statement of the main theorem.  In Section \ref{crysec}, we give the necessary background information on crystals.  Section \ref{pmsec} gives a brief review of $\pm$-sequences, which are at the heart of our crystal constructions.  Section \ref{sl3sec} contains an abstract construction of crystals which turn out to be crystals of highest-weight $\fsl_3$-modules; these crystals are used in Section \ref{proofsec} to complete the proof of the main theorem.  Finally in Section \ref{diffsec} we show that the various models we have created for $B(\La_0)$ for a given value of $n$ are distinct, so that we really do get an uncountable family.

\begin{acks}
The author's thanks are due to Chris Berg, who made an early version of his paper \cite{berg} available.  The author enjoyed fruitful discussions with Chris (as well as Monica Vazirani and Brant Jones) at MSRI Berkeley in March 2008, during the concurrent MSRI programmes `Combinatorial representation theory' and `Representation theory of finite groups and related topics'.  The author would like to thank the organisers of these programmes, as well as MSRI for some financial support.

This work was begun while the author was a visiting Postdoctoral Fellow at Massachusetts Institute of Technology, with the aid of a Research Fellowship from the Royal Commission for the Exhibition of 1851.  The author is very grateful to M.I.T.\ for its hospitality, and the 1851 Commission for its generous support.
\end{acks}

\sect{The main theorem}

In this section we state our main theorem.  In order to get to the result as quickly as possible, we give only minimal definitions here; we assume the reader has some familiarity with representations of Kac--Moody algebras and quantum groups.  Further definitions are given in later sections.

\subs{The algebra $\hsl_n$}

In this paper we fix an integer $n\gs3$, and consider the Kac--Moody algebra $\hsl_n$.  This has simple roots indexed by the set $I=\znz$, and the Cartan matrix is given by
\[a_{ij} = 2\delta_{ij}-\delta_{i(j+1)}-\delta_{i(j-1)}\]
for $i,j\in I$.  We may abuse notation and label elements of $\znz$ by the integers $0,\dots,n-1$; of course, the subscripts in the above formula should then be read modulo $n$.

The usual realisation of this Cartan matrix is given as follows.  The Cartan subalgebra $\fkh$ is an $(n+1)$-dimensional vector space over $\bbc$ with basis
\[\{h_i\mid i\in\znz\}\cup\{D\}.\]
The dual space $\fkh^\ast$ has basis
\[\{\Lambda_i\mid i\in\znz\}\cup\{\delta\},\]
and the weight lattice $P$ is the $\bbz$-span of this basis.  The two bases are paired via
\[\lan h_i,\Lambda_j\ran = \delta_{ij},\qquad \lan D,\delta\ran=1,\qquad \lan h_i,\delta\ran=\lan D,\Lambda_j\ran=0,\]
and the simple roots $\alpha_i$ ($i\in\znz$) are given by
\[\alpha_i = 2\Lambda_i-\Lambda_{i-1}-\Lambda_{i+1}+\delta_{i0}\delta.\]
In this paper, we shall be concerned with the highest-weight module $V(\Lambda_0)$ for the quantum enveloping algebra $U_q(\hsl_n)$; as well as being a prototype for highest-weight modules of affine Kac--Moody algebras, this has important applications in the representation theory of Iwahori--Hecke algebras; for example, see \cite{groj}.

\subs{Partitions}

As usual, a \emph{partition} is a weakly decreasing sequence $\la=(\la_1,\la_2,\dots)$ of non-negative integers such that the sum $\la_1+\la_2+\dots$ is finite.  Partitions are often written with equal parts grouped together and zeroes omitted, and the partition $(0,0,\dots)$ is written as $\varnothing$.  A partition $\la$ is often identified with its \emph{Young diagram}, which is the set
\[[\la]=\left\{\left.(a,c)\in\bbn^2\ \right|\ c\ls\la_a\right\}.\]
We adopt the English convention for drawing Young diagrams, in which $a$ increases down the page and $c$ increases from left to right.

We refer to elements of $\bbn^2$ as \emph{nodes}, and elements of $[\la]$ as nodes of $\la$.  A node $(a,c)$ of $\la$ is \emph{removable} if $[\la]\setminus\{(a,c)\}$ is the Young diagram of a partition (i.e.\ if $c=\la_a>\la_{a+1}$), while a node $(a,c)$ not in $\la$ is an \emph{addable node} of $\la$ if $[\la]\cup\{(a,c)\}$ is the Young diagram of a partition.

If $\la$ is a partition, the \emph{conjugate partition} $\la'$ is defined by
\[\la'_a = \big|\left\{c\ \left|\ \la_c\gs a\right.\right\}\big|.\]

Now suppose we have fixed $n\gs2$.  We define the \emph{residue} of a node $(a,c)$ to be $c-a+n\bbz$.  If $(a,c)$ has residue $i\in\znz$, we refer to it as an \emph{$i$-node}.

\subs{Partition models for $B(\Lambda_0)$}\label{models}

Now we describe the two usual partition models for the basic crystal $B(\La_0)$ for $U_q(\hsl_n)$; see Section \ref{crysec} for an introduction to crystals.

Say that a partition $\la$ is \emph{$n$-restricted} if $\la_a-\la_{a+1}<n$ for all $a$.  We write $\rest n$ for the set of $n$-restricted partitions, and we define crystal operators $\ee i,\ff i$ for $i\in\znz$ as follows.  Given $\la\in\rest n$, let $(a_1,c_1),\dots,(a_r,c_r)$ be the list of all addable and removable $i$-nodes of $\la$, ordered so that $c_1<\dots<c_r$.  Now write a sequence of signs $\pi=\pi_1\dots\pi_r$, where $\pi_j$ equals $+$ if $(a_j,c_j)$ is an addable node, and $-$ otherwise.

If there is no good position in this sequence (see \S\ref{pmsec} for the definition of a good position), then we set $\ee i\la=0$.  Otherwise, if position $j$ is the good position, we define $\ee i\la$ to be the partition obtained from $\la$ by removing the node $(a_j,c_j)$.  Similarly, if there is no cogood position in $\pi$, then we define $\ff i\la=0$, and otherwise we define $\ff i\la$ to be the partition obtained by adding the node corresponding to the cogood position in $\pi$.

We define a function $\wt$ on partitions by
\[\wt(\la) = \La_0-\sum_{i\in\znz}c_i\alpha_i,\]
where $c_i$ is the number of $i$-nodes of $\la$.  Then we have the following.

\begin{thm}\thmcite{mm}{Theorem 4.7}\label{restcrystal}
The operators $\ee i,\ff i$ and the weight function $\wt$ endow $\rest n$ with the structure of a crystal for $\hsl_n$.  This crystal is isomorphic to $B(\Lambda_0)$.
\end{thm}

Now we describe the other well-known realisation of the basic $U_q(\hsl_n)$-crystal in terms of partitions.  Say that a partition $\la$ is \emph{$n$-regular} if there is no $a$ for which $\la_a=\la_{a+n-1}>0$; in other words, $\la$ is $n$-regular if $\la'$ is $n$-restricted.  We may define crystal operators $\ee i,\ff i$ on the set of $n$-regular partitions; this is done exactly as for $n$-restricted partitions above, except that we replace $c_1<\dots<c_r$ above with $c_1>\dots>c_r$.  Then we get a counterpart to Theorem \ref{restcrystal}.  An explicit isomorphism between these two crystals may be given in terms of the Mullineux map \cite{mull}.

\vspace{\topsep}
The object of the present paper is to describe, for each $n\gs3$, an uncountable family of realisations of the crystal $B(\La_0)$ in terms of partitions.  The $n$-regular and $n$-restricted versions described above arise as special cases, as does Berg's `ladder crystal' \cite{berg}, which was the starting point for the present work.  Since the first version of this paper was written, Tingley \cite{ting} has shown that another model from the literature -- Nakajima's monomial crystal -- occurs as a special case too.

We now fix $n\gs3$, and describe our crystals.  The rest of the paper is devoted to proving that these crystals are all isomorphic to the crystal $B(\La_0)$.

Each of our crystals will be indexed by an \emph{arm sequence}, which is defined to be a sequence $A=A_1,A_2,\dots$ of integers such that
\begin{itemize}
\item
$t-1\ls A_t\ls (n-1)t$ for all $t\gs1$, and
\item
$A_{t+u}\in\{A_t+A_u,A_t+A_u+1\}$ for all $t,u\gs1$.
\end{itemize}

In order to give a description of the particular partitions which will label the vertices of our crystal, we need to discuss hooks.  Recall that if $\la$ is any partition and $(a,c)$ is a node of $\la$, then the \emph{$(a,c)$-hook} of $\la$ is the set of nodes of $\la$ directly to the right of or directly below $(a,c)$, including $(a,c)$ itself.  The \emph{$(a,c)$-hook length} of $\la$ is the number of nodes in the $(a,c)$-hook, i.e.\ the number
\[\la_a-c+\la'_c-a+1,\]
while the \emph{$(a,c)$-arm length} is just the number of nodes strictly to the right of $(a,c)$, i.e.\ $\la_a-c$.

Now suppose $A$ is an arm sequence.  If $\la$ is a partition, we say that $\la$ is \emph{$A$-regular} if there is no node $(a,c)$ of $\la$ such that the $(a,c)$-hook length of $\la$ equals $nt$ and the $(a,c)$-arm length equals $A_t$, for any $t$.

We refer to a hook with length $nt$ and arm length $A_t$ as an \emph{illegal hook}; so an $A$-regular partition is one with no illegal hooks.

\begin{egs}\indent
\vspace{-\topsep}
\begin{enumerate}
\item
Suppose $A_t=(n-1)t$ for all $t$.  Then a partition $\la$ is $A$-regular if and only if it is $n$-restricted.  To see this, suppose first that $\la$ is not $n$-restricted.  Then we have $\la_a-\la_{a+1}\gs n$ for some $a$, so the $(a,\la_a-n+1)$-hook of $\la$ has length $n$ and arm length $n-1$.  Conversely, suppose that for some node $(a,c)$ of $\la$ and $t\gs1$ the $(a,c)$-hook length of $\la$ equals $nt$ while the $(a,c)$-arm length equals $(n-1)t$.  Then we have $\la_a=c+(n-1)t$ and $\la_{a+t}<c$, so $\la_a-\la_{a+t}>(n-1)t$, and hence we must have $\la_b-\la_{b+1}\gs n$ for some $b\in\{a,a+1,\dots,a+t-1\}$; so $\la$ is not $n$-restricted.
\item
Suppose $A_t=t$ for each $t$; then a partition is $A$-regular if and only if it does not possess a hook whose length is exactly $n$ times its arm length.  Hence by \cite[Theorem 13.2.3]{berg} the $A$-regular partitions are precisely the partitions appearing in Berg's ladder crystal.  This characterisation of these partitions in terms of illegal hooks was the starting point for the present work.
\end{enumerate}
\end{egs}

Write $\regg A$ for the set of $A$-regular partitions; this will be the underlying set of our crystal.  In order to define the crystal operators $\ee i,\ff i$ on $\regg A$, we need to introduce a total order on the set of all nodes of a given residue. If $(a,c)$ and $(b,d)$ are distinct nodes with the same residue, then the \emph{axial distance} $b-a+c-d$ equals $nt$ for some integer $t$; by interchanging $(a,c)$ and $(b,d)$ if necessary, we suppose $t\gs0$.  Now we set $(a,c)\cg(b,d)$ if $c-d>A_t$, and $(b,d)\cg(a,c)$ otherwise; for this purpose, we read $A_0$ as $0$.

It easy to check, using the definition of an arm sequence, that $\cg$ is transitive, and is therefore a total order on the set of nodes of a given residue.  Now we can define our crystal operators.  Suppose $\la\in\regg A$ and $i\in\znz$.  Let $(a_1,c_1),\dots,(a_r,c_r)$ be the list of addable and removable $i$-nodes of $\la$, ordered so that $(a_1,c_1)\cg\cdots\cg(a_r,c_r)$.  Define the sequence of signs $\pi=\pi_1\dots\pi_r$ by putting $\pi_j=+$ if $(a_j,c_j)$ is an addable node, and $\pi_j=-$ otherwise.  If $\pi$ has no good position, then set $\ee i\la=0$.  Otherwise, let $j$ be the good position, and say that $(a_j,c_j)$ is the $i$-good node for $\la$; define $\ee i\la$ by removing this node from $\la$.  Similarly, if there is no cogood position in $\pi$, then set $\ff i\la=0$, and otherwise let $k$ be the cogood position and say that $(a_k,c_k)$ is the $i$-cogood node for $\la$; define $\ff i\la$ by adding this node to $\la$.

Now we can state our main theorem.

\begin{thm}\label{main}
Suppose $n\gs3$ and $A$ is an arm sequence.  Then the crystal operators $\ee i,\ff i$ and the weight function $\wt$ endow $\regg A$ with the structure of an $\hsl_n$-crystal.  This crystal is isomorphic to $B(\La_0)$.
\end{thm}

\begin{eg}
Suppose $n=3$, and $A$ satisfies $A_1=1$, $A_2=3$ and $A_3=4$.  Part of the crystal $\regg A$ is shown in Figure 1.
\begin{figure}
\vspace{-10pt}
\Yboxdim{4.4pt}

\makeatletter
\def\dggeometry{%
   \dg@ZTEMP=\dg@XGRID \multiply\dg@ZTEMP\tw@
   \ifnum\dg@YGRID=\z@ \dg@ZTEMP=\tw@
   \else \divide\dg@ZTEMP\dg@YGRID \fi
   \ifnum\dg@ZTEMP>\f@ur \dg@ZTEMP=\f@ur \fi
   \ifnum\dg@ZTEMP<\@ne \dg@ZTEMP=\@ne \fi
   \unitlength=1sp\relax
\dg@XGRID=4
\dg@YGRID=3
\unitlength=840sp
}
\makeatother

\dgLABELOFFSET=1pt
\dgARROWLENGTH=10pt
\dgHORIZPAD=2pt
\dgVERTPAD=4pt
\[\begin{diagram}
\ny{10}9\\
\ny98\ane2\ny1{8,1}\\
\ny87\ane1\are2\ny1{7,1}\ane1\are0\ny1{7,2}\\
\ny76\ane0\are2\ny1{6,1}\ane0\arrow{se,t,1}1\ny1{6,2}\arrow{se,t,1}1\ny1{7,1,1}\\
\ny65\ane2\ny1{4,2}\are1\ny1{5,2}\arrow{ne,t,1}2\ase1\ny1{6,1,1}\arrow{ne,t,1}0\ny1{6,3}\\
\ny54\ane1\ase2\ny4{5,3}\ane2\are1\ny1{5,3,1}\\
\ny43\ane0\ase2\ny2{4,1}\are1\arrow{nne,t}0\ny1{5,1}\arrow{nne,t}0\are1\ny1{5,1,1}\arrow{nne,t}2\ase0\ny2{6,1,1,1}\\
\ny32\ane2\ny2{3,1}\ane0\ase1\ny4{5,1,1,1}\ane2\are0\ny1{5,2,1,1}\\
\node{\varnothing}\are0\ny11\ane1\ase2\ny4{3,1,1}\are0\ny1{3,1,1,1}\are0\ny1{4,1,1,1}\ane1\are0\ase2\ny1{4,2,1,1}\ane1\ase2\\
\ny3{1,1}\ase1\ny2{2,1,1}\ase0\ane2\ny4{4,1,1,1,1}\ase1\are0\ny1{4,2,1,1,1}\\
\ny4{1,1,1}\ase0\ane1\ny2{2,1,1,1}\are2\arrow{sse,t}0\ny1{2,1,1,1,1}\arrow{sse,t}0\are2\ny1{3,1,1,1,1}\arrow{sse,t}1\ane0\ny2{4,1,1,1,1,1}\\
\ny5{1,1,1,1}\ane1\ase2\ny4{2,2,2,1,1}\ase1\are2\ny1{3,2,2,1,1}\\
\ny6{1,1,1,1,1}\ase1\ny1{2,2,1,1}\are2\ny1{2,2,1,1,1}\arrow{se,t,1}1\ane2\ny1{3,1,1,1,1,1}\arrow{se,t,1}0\ny1{2,2,2,1,1,1}\\
\ny7{1,1,1,1,1,1}\ase0\are1\ny1{2,1,1,1,1,1}\ase0\arrow{ne,t,1}2\ny1{2,2,1,1,1,1}\arrow{ne,t,1}2\ny1{3,1,1,1,1,1,1}\\
\ny8{1,1,1,1,1,1,1}\ase2\are1\ny1{2,1,1,1,1,1,1}\ase2\are0\ny1{2,2,1,1,1,1,1}\\
\ny9{1,1,1,1,1,1,1,1}\ase1\ny1{2,1,1,1,1,1,1,1}\\
\ny{10}{1,1,1,1,1,1,1,1,1}
\end{diagram}\]
\vspace{-20pt}
\begin{center}Figure 1\end{center}
\end{figure}
\end{eg}

\sect{Background on crystals}\label{crysec}

In this section, we give the definitions and basic results we shall need concerning crystals.  An indispensable introduction to this subject is Kashiwara's book \cite{kash}.

\subs{Crystals}\label{crydef}

Suppose $I$ is a finite set.  We define an \emph{$I$-crystal} to be a finite set $B$, equipped with functions
\[\ee i,\ff i:B\lra B\sqcup\{0\}\]
for each $i\in I$, satisfying the following axioms.
\renewcommand{\theenumi}{C\arabic{enumi}}
\begin{enumerate}
\item
If $b,b'\in B$, then
\[b=\ee ib'\qquad\Longleftrightarrow\qquad b'=\ff ib.\]
\item
If we set
\[\ep_i(b)=\max\big\{\ep\ \big|\ \ee i^\ep(b)\neq0\big\},\qquad\phi_i(b)=\max\big\{\phi\ \big|\ \ff i^\phi(b)\neq0\big\},\]
then $\ep_i(b),\phi_i(b)$ are finite.
\end{enumerate}
\renewcommand{\theenumi}{\arabic{enumi}}
Here $0$ is what Kashiwara calls a `ghost element'.  In other words, $\ee i$ and $\ff i$ may be regarded as mutually inverse injective partial functions from $B$ to itself.  To an $I$-crystal is associated a \emph{crystal graph}: this has vertex set $B$, with an arrow labelled with $i\in I$ from $b$ to $\ff ib$ whenever $\ff ib\neq0$.  An $I$-crystal is often identified with its crystal graph, and accordingly one may speak of a crystal being connected, or talk of the connected components of a crystal.  Of course, each connected component is itself a crystal.

Now let $\fkg$ be a symmetrisable Kac--Moody algebra.  This is defined by a Cartan matrix $A=(a_{ij})_{i,j\in I}$,
where $I$ is an indexing set, which we assume to be finite.  For each $i\in I$, we let $\alpha_i$ denote the corresponding simple root, and $h_i$ the simple coroot.  $P$ denotes the lattice of integral weights.  A $\fkg$-\emph{crystal} is defined to be an $I$-crystal, equipped with a function
\[\wt:B\lra P\]
such that the following additional axioms are satisfied for each $i\in I$.
\renewcommand{\theenumi}{C\arabic{enumi}}
\begin{enumerate}
\setcounter{enumi}2
\item
If $b\in B$ and $\ee ib\neq0$, then
\[\wt(\ee ib) = \wt(b)+\alpha_i.\]
\item\label{c4}
For each $b\in B$,
\[\lan h_i,\wt(b)\ran = \phi_i(b)-\ep_i(b).\]
\end{enumerate}
\renewcommand{\theenumi}{\arabic{enumi}}

\begin{rmks}\indent
\vspace{-\topsep}
\begin{enumerate}
\item
In fact, the definition we have given above is the definition of a \emph{semi-normal} crystal; in general, a more liberal definition of $\ep_i,\phi_i$ is permitted.  Since in this paper we shall only be concerned with semi-normal crystals, we use the term `crystal' as defined above.
\item
Note that if the Cartan matrix is non-singular (in particular, if $\fkg$ is of finite type), then by axiom \ref{c4} the functions $\ee i,\ff i$ determine the weight function $\wt$.  Hence in this case we may specify a crystal simply by giving the functions $\ee i,\ff i$ for each $i$; the assertion that we have a crystal is then the assertion that these functions, together with the implied weight function, satisfy the axioms.
\end{enumerate}
\end{rmks}

\subsubsection*{Subcrystals}

Given an $I$-crystal $B$ and a subset $J$ of $I$, one obtains a $J$-crystal $B_J$ simply by forgetting the functions $\ee i,\ff i$ for $i\in I\setminus J$; we call this a \emph{subcrystal} of $B$.  If $B$ is a $\fkg$-crystal with $\fkg$ having indexing set $I$, then $B_J$ is a $\fkg_J$-crystal, where $\fkg_J$ is the subalgebra of $\fkg$ corresponding to $J$; the weight function on $B_J$ is induced from the weight function on $B$.

We shall be particularly interested in the case where $|J|=2$; in this case we refer to $B_J$ as a \emph{rank $2$ subcrystal} of $B$.

\subsubsection*{Tensor products}

If $B,B'$ are $I$-crystals, the \emph{tensor product} $B\otimes B'$ is an $I$-crystal defined as follows.  The underlying set is the Cartesian product $B\times B'$, whose elements we write in the form $b\otimes b'$, for $b\in B$ and $b'\in B'$.  For $i\in I$ the functions $\ee i,\ff i$ are given by
\[\ee i(b\otimes b')=\begin{cases}
(\ee ib)\otimes b' & (\phi(b)\gs\ep(b'))\\
b\otimes (\ee ib') & (\phi(b)<\ep(b')),
\end{cases}\qquad
\ff i(b\otimes b')=\begin{cases}
(\ff ib)\otimes b' & (\phi(b)>\ep(b'))\\
b\otimes (\ff ib') & (\phi(b)\ls\ep(b'));
\end{cases}\]
in these definitions, $b\otimes0$ and $0\otimes b'$ should both be taken to equal $0$.

If $B,B'$ are $\fkg$-crystals, then $B\otimes B'$ is also a $\fkg$-crystal, with weight function
\[\wt(b\otimes b') = \wt(b)+\wt(b').\]
The point of this definition is that if $B,B'$ are crystals arising from crystal bases of $U_q(\fkg)$-modules $V,V'$, then $B\otimes B'$ is the crystal given by the associated crystal basis of $V\otimes V'$ (which is a $U_q(\fkg)$-module once an appropriate coproduct on $U_q(\fkg)$ is chosen).

\subsection{Crystals of highest-weight modules}

Of particular interest in the study of $\fkg$-crystals are highest-weight crystals.  For each dominant integral weight $\La$ of $\fkg$, Kashiwara \cite{ka1} proved that the irreducible highest-weight module $V(\La)$ has an essentially unique crystal basis; we write $B(\La)$ for the corresponding $\fkg$-crystal.  Following Danilov et al.\ \cite{dkk}, we say that a $\fkg$-crystal is \emph{regular} if it is isomorphic to $B(\La)$ for some dominant integral weight $\La$.

The aim of this paper is to give a family of constructions of one particular regular crystal $B(\La_0)$.  Our approach to this is to define an abstract crystal, and then to prove that it is regular and read off the highest weight.  In order to do this, we shall need to use the following results.

\begin{propn}\label{tenshw}
Suppose $B,B'$ are regular $\fkg$-crystals.  Then every connected component of $B\otimes B'$ is a regular crystal.
\end{propn}

\begin{pf}
This follows from the fact that a tensor product of two irreducible highest-weight modules decomposes as a direct sum of irreducible highest-weight modules.
\end{pf}

Now say that an element $b$ of an $I$-crystal $B$ is a \emph{source} if $\ee ib=0$ for all $i\in I$.

\begin{thm}\label{rank2}\thmcite{kkmmnn}{Proposition 2.4.4}
Suppose $B$ is a $\fkg$-crystal.  Then $B$ is regular if and only if $B$ has a unique source and every connected component of every rank $2$ subcrystal of $B$ is regular.
\end{thm}

Once we know that a crystal is regular, it is straightforward to find its highest weight; this is just the weight of the unique source.

In view of Theorem \ref{rank2}, it will be useful to have a characterisation of regular crystals in the case where $\fkg$ has rank $2$.  In fact, the only cases which will interest us are those which arise as rank $2$ subalgebras of $\hsl_n$ for $n\gs3$, namely $\fsl_2\oplus\fsl_2$ and $\fsl_3$.  These cases were studied by Stembridge \cite{stem}; he gave a list of `local' axioms which characterise highest-weight crystals in these cases.  In fact, his results are stronger: starting from an $I$-crystal $B$ where $|I|=2$ (and not \textit{a priori} a crystal for any particular $\fkg$) his axioms characterise whether $B$ is a regular crystal for $\fsl_2\oplus\fsl_2$ or $\fsl_3$ (with the implied weight function -- see Remark 2 in \S\ref{crydef}).  Using Littelmann's path model, Stembridge goes on to give a new proof of Theorem \ref{rank2} in the case where $\fkg$ is simply-laced.

\subsubsection*{Highest-weight crystals for $\fsl_2\oplus\fsl_2$}

Suppose $|I|=2$, say $I=\{i,j\}$, and that the Cartan matrix $A$ is given by
\[a_{ii}=a_{jj}=2,\qquad a_{ij}=a_{ji}=0;\]
then the associated Kac--Moody algebra is isomorphic to $\fsl_2\oplus\fsl_2$, and we identify it with the latter.  The following result is a special case of the results in \cite{stem}.

\begin{propn}\label{sl2}
Suppose $B$ is a connected $I$-crystal, where $I=\{i,j\}$.  Then $B$ is a regular $\fsl_2\oplus\fsl_2$-crystal if and only if $\ee i,\ff i$ both commute with $\ee j,\ff j$.
\end{propn}

What this says is that $B$ is a regular $\fsl_2\oplus\fsl_2$-crystal if and only if its crystal graph is `rectangular', as in the following diagram:
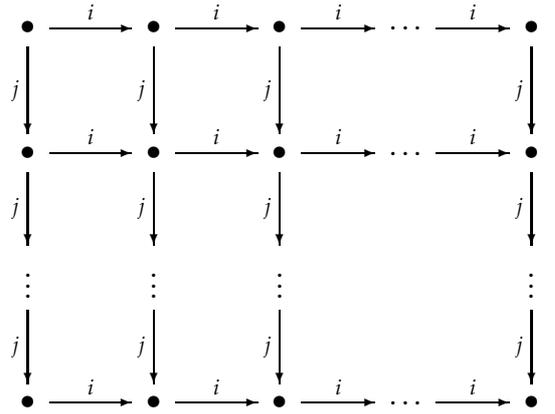
\[\begin{diagram}
\node{\bullet}\arrow{e,t}i\arrow{s,l}j\node{\bullet}\arrow{e,t}i\arrow{s,l}j\node{\bullet}\arrow{e,t}i\arrow{s,l}j\node{\dots}\arrow{e,t}i\node{\bullet}\arrow{s,l}j\\
\node{\bullet}\arrow{e,t}i\arrow{s,l}j\node{\bullet}\arrow{e,t}i\arrow{s,l}j\node{\bullet}\arrow{e,t}i\arrow{s,l}j\node{\dots}\arrow{e,t}i\node{\bullet}\arrow{s,l}j\\
\node{\vdots}\arrow{s,t}j\node{\vdots}\arrow{s,t}j\node{\vdots}\arrow{s,t}j\node[2]{\vdots}\arrow{s,l}j\\
\node{\bullet}\arrow{e,t}i\node{\bullet}\arrow{e,t}i\node{\bullet}\arrow{e,t}i\node{\dots}\arrow{e,t}i\node{\bullet}
\end{diagram}.\]

\subsubsection*{Highest-weight crystals for $\fsl_3$}

If $I=\{i,j\}$, with
\[a_{ii}=a_{jj}=2,\qquad a_{ij}=a_{ji}=-1,\]
then we shall identify $\fkg$ with $\fsl_3$.  Stembridge's axioms to characterise regular $\fsl_3$-crystals are complicated, and it seems to be difficult to verify these directly for the $\fsl_3$-crystals in this paper.  Accordingly, we shall use a different version, due to Danilov, Karzanov and Koshevoy.

Suppose we are given an $I$-crystal $B$, and define the functions $\ep_i,\phi_i,\ep_j,\phi_j$ on $B$ as in \S\ref{crydef}.  Now we may state axioms A2 and A3 from \cite{dkk} as follows (their axiom A1 is implicit in our definition of a crystal).

\renewcommand{\theenumi}{A\arabic{enumi}}
\begin{enumerate}
\setcounter{enumi}1
\item
\begin{enumerate}
\item
If $b\in B$ and $\ee i(b)\neq0$, then either
\begin{alignat*}2
\ep_j(\ee ib) &= \ep_j(b)+1,&\qquad\phi_j(\ee ib)&=\phi_j(b)\\
\intertext{or}
\ep_j(\ee ib) &= \ep_j(b),&\qquad\phi_j(\ee ib)&=\phi_j(b)-1.
\end{alignat*}
\item
If in addition $\ff ib\neq0$, then either
\[\ep_j(\ee ib)>\ep_j(b)\qquad\text{or}\qquad\phi_j(\ff ib)>\phi_j(b).\]
\end{enumerate}
The same statements hold with $i$ and $j$ interchanged.
\item
\begin{enumerate}
\item
If $b\in B$, $\ee ib\neq0$, $\ee jb\neq0$ and $\ep_j(\ee ib)=\ep_j(b)$, then
\[\ee i\ee jb=\ee j\ee ib\neq0\qquad\text{and}\qquad\ep_i(\ee jb)>\ep_i(b).\]
\item
If $b\in B$, $\ff ib\neq0$, $\ff jb\neq0$ and $\phi_j(\ff ib)=\phi_j(b)$, then
\[\ff i\ff jb=\ff j\ff ib\neq0\qquad\text{and}\qquad\phi_i(\ff jb)>\phi_i(b).\]
\end{enumerate}
The same statements hold with $i$ and $j$ interchanged.
\end{enumerate}
\renewcommand{\theenumi}{\arabic{enumi}}

There is also an axiom A4 in \cite{dkk}, but we shall not need this, since we use the following theorem.

\begin{propn}\label{aaas}\thmcite{dkk}{Proposition 5.3}
Suppose $B$ is an $\{i,j\}$-crystal satisfying axioms A2 and A3, and suppose that $B$ has a unique source.  Then $B$ is a regular $\fsl_3$-crystal.
\end{propn}

\sect{$\pm$-sequences}\label{pmsec}

In this section, we give some basic properties of $\pm$-sequences; these are essential to the definitions of our crystals.

We define a \emph{$\pm$-sequence} to be a finite word $\pi=\pi_1\dots \pi_m$ from the alphabet $\{+,0,-\}$.  We shall perform arithmetic with these signs, interpreting $+$ as $+1$ and $-$ as $-1$.

If $\pi_1\dots \pi_m$ is a $\pm$-sequence, we define
\[h_i(\pi) = \pi_i+\dots+\pi_m\]
for $i=1,\dots,m+1$.  We define
\[\phi(\pi) = \max\left\{\left.h_i(\pi)\ \rule{0cm}{9pt}\right|\ 1\ls i\ls m+1\right\},\]
and if $\phi(\pi)>0$ we define the \emph{cogood} position in $\pi$ to be
\[\max\left\{i\ \left|\ h_i(\pi)=\phi(\pi)\right.\right\}.\]

Similarly, we define
\[g_i(\pi) = -(\pi_1+\dots+\pi_i)\]
for $i=0,\dots,m$, and we let
\[\ep(\pi) = \max\left\{\left.g_i(\pi)\ \right|\ 0\ls i\ls m\right\};\]
if $\ep(\pi)>0$ we define the \emph{good} position in $\pi$ to be
\[\min\left\{i\ \left|\ g_i(\pi)=\ep(\pi)\right.\right\}.\]

\begin{eg}
Suppose
\begin{alignat*}2
\pi&=&&\makebox[\spacer]{$-$}\makebox[\spacer]{$+$}\makebox[\spacer]{$0$}\makebox[\spacer]{$-$}\makebox[\spacer]{$-$}\makebox[\spacer]{$0$}\makebox[\spacer]{$0$}\makebox[\spacer]{$+$}\makebox[\spacer]{$+$}\makebox[\spacer]{$0$}\makebox[\spacer]{$-$}\makebox[\spacer]{$0$}\makebox[\spacer]{$-$}\makebox[\spacer]{$+$}\makebox[\spacer]{$+$}\makebox[\spacer]{$0$}\makebox[\spacer]{$-$}.\\
\intertext{Then the values of the functions $h$ and $g$ are as follows.}
h&=&&\makebox[\spacer]{$-\negthinspace1$}\makebox[\spacer]{$0$}\makebox[\spacer]{$-\negthinspace1$}\makebox[\spacer]{$-\negthinspace1$}\makebox[\spacer]{$0$}\makebox[\spacer]{$1$}\makebox[\spacer]{$1$}\makebox[\spacer]{$1$}\makebox[\spacer]{$0$}\makebox[\spacer]{$-\negthinspace1$}\makebox[\spacer]{$-\negthinspace1$}\makebox[\spacer]{$0$}\makebox[\spacer]{$0$}\makebox[\spacer]{$1$}\makebox[\spacer]{$0$}\makebox[\spacer]{$-\negthinspace1$}\makebox[\spacer]{$-\negthinspace1$}\makebox[\spacer]{$0$}\\
g&=&\makebox[\spacer]{$0$}&\makebox[\spacer]{$1$}\makebox[\spacer]{$0$}\makebox[\spacer]{$0$}\makebox[\spacer]{$1$}\makebox[\spacer]{$2$}\makebox[\spacer]{$2$}\makebox[\spacer]{$2$}\makebox[\spacer]{$1$}\makebox[\spacer]{$0$}\makebox[\spacer]{$0$}\makebox[\spacer]{$1$}\makebox[\spacer]{$1$}\makebox[\spacer]{$2$}\makebox[\spacer]{$1$}\makebox[\spacer]{$0$}\makebox[\spacer]{$0$}\makebox[\spacer]{$1$}.
\end{alignat*}
Hence $\ep(\pi)=2$ and $\phi(\pi)=1$, with the good and cogood positions being $5$ and $14$ respectively.
\end{eg}

\begin{rmks}\indent
\begin{enumerate}
\vspace{-\topsep}
\item
The reader will observe that the zeroes in a $\pm$-sequence $\pi$ have no effect on $\phi(\pi)$ or $\ep(\pi)$.  That is, if we define a new sequence $\pi'$ by removing some or all of the zeroes, then we have $\phi(\pi)=\phi(\pi')$ and $\ep(\pi)=\ep(\pi')$, and the good and cogood positions in $\pi'$ correspond to those in $\pi$.  It is slightly more useful to us in this paper to allow zeroes.
\item
The reader may prefer an alternative way to work out $\ep(\pi)$ and $\phi(\pi)$ and the good and cogood positions.  If $\pi$ is a $\pm$-sequence, then we define the \emph{reduction} of $\pi$ by repeatedly replacing any segment of the form $+00\dots0-$ with a segment $00\dots0$ of the same length until there are no more such segments (i.e.\ there is no $+$ preceding a $-$).  Then $\ep(\pi)$ is the number of $-$ signs in the reduction of $\pi$, and the good position is the position of the last such sign; $\phi(\pi)$ is the number of $+$ signs in the reduction of $\pi$, and the cogood position is the position of the first $+$.  For instance, the reduction of the sequence $\pi$ given in the above example is the following:
\[\makebox[\spacer]{$-$}\makebox[\spacer]{$0$}\makebox[\spacer]{$0$}\makebox[\spacer]{$0$}\makebox[\spacer]{$-$}\makebox[\spacer]{$0$}\makebox[\spacer]{$0$}\makebox[\spacer]{$0$}\makebox[\spacer]{$0$}\makebox[\spacer]{$0$}\makebox[\spacer]{$0$}\makebox[\spacer]{$0$}\makebox[\spacer]{$0$}\makebox[\spacer]{$+$}\makebox[\spacer]{$0$}\makebox[\spacer]{$0$}\makebox[\spacer]{$0$}.\]
\end{enumerate}
\end{rmks}

We now collect some basic properties of $\pm$-sequences.  We omit the proofs in the hope that the reader will find it more instructive to prove the results himself.

\begin{lemma}\label{pmutil}
Suppose $\pi=\pi_1\dots \pi_m$ is a $\pm$-sequence.
\begin{enumerate}
\item\label{pmu1}
$\phi(\pi)-\ep(\pi)=\pi_1+\dots+\pi_m$.
\item\label{pmu2}
If $\ep(\pi)>0$ and $i$ is the good position in $\pi$, then $\pi_i=-$.
\item\label{pmu3}
If $\phi(\pi)>0$ and $j$ is the cogood position in $\pi$, then $\pi_j=+$.
\item\label{pmu4}
If $\ep(\pi),\phi(\pi)>0$ and $i,j$ are the good and cogood positions respectively, then $i<j$.
\item\label{pmu5}
If $i$ is the good position in $\pi$ and $j>i$ with $\pi_j=-$, then there exists $k$ such that $i<k<j$ and $\pi_k=+$.
\item\label{pmu6}
If $j$ is the cogood position in $\pi$ and $i<j$ with $\pi_i=+$, then there exists $k$ such that $i<k<j$ and $\pi_k=-$.
\end{enumerate}
\end{lemma}

Now we consider the relationship between two $\pm$-sequences.

\begin{lemma}\label{goodcogood}
Suppose $\pi=\pi_1\dots \pi_m$ and $\rho=\rho_1\dots \rho_m$ are two $\pm$-sequences, and $i\in\{1,\dots,m\}$ is such that
\[\pi_i=-,\qquad \rho_i=+,\qquad \pi_j=\rho_j\ \text{ for all }j\neq i.\]
Then position $i$ is the good position in $\pi$ if and only if it is the cogood position in $\rho$.  If this is the case, then
\[\ep(\rho)=\ep(\pi)-1,\qquad \phi(\rho)=\phi(\pi)+1.\]
\end{lemma}

\begin{lemma}\label{downbyone}
Suppose $\pi=\pi_1\dots \pi_m$ and $\rho=\rho_1\dots \rho_m$ are two $\pm$-sequences, and $i\in\{1,\dots,m\}$ is such that
\[\rho_i=\pi_i-1,\qquad \pi_j=\rho_j\ \text{ for all }j\neq i.\]
Then we have the following.
\begin{enumerate}
\item
If $\phi(\pi)>0$ and the cogood position in $\pi$ is at or before position $i$, then
\[\epsilon(\rho)=\epsilon(\pi),\qquad \phi(\rho)=\phi(\pi)-1.\]
  If in addition $\epsilon(\pi)>0$, then $\pi$ and $\rho$ have the same good position.
\item
If $\phi(\pi)=0$ or the cogood position in $\pi$ is strictly after position $i$, then
\[\epsilon(\rho)=\epsilon(\pi)+1,\qquad \phi(\rho)=\phi(\pi).\]
\item
If $\epsilon(\pi)>0$ and the good position in $\pi$ is after position $i$, then $\pi$ and $\rho$ have the same good position.
\end{enumerate}
\end{lemma}

Finally, we need to consider the concatenation of two $\pm$-sequences.

\begin{lemma}\label{concat}
Suppose $\pi=\pi_1\dots \pi_l$ and $\rho=\rho_1\dots \rho_m$ are two $\pm$-sequences, and let
\[\pi\ast\rho = \pi_1\dots \pi_l\rho_1\dots \rho_m\]
denote their concatenation.
\begin{itemize}
\item
Suppose $\phi(\pi)\gs \ep(\rho)$.  Then $\ep(\pi\ast\rho)=\ep(\pi)$, and if this is positive, then the good position in $\pi\ast\rho$ coincides with the good position in $\pi$.
\item
Suppose $\phi(\pi)<\ep(\rho)$.  Then $\ep(\pi\ast\rho)=\ep(\pi)+\ep(\rho)-\phi(\pi)$, and the good position in $\pi\ast\rho$ corresponds to the good position in $\rho$.
\item
Suppose $\phi(\pi)\ls \ep(\rho)$.  Then $\phi(\pi\ast\rho)=\phi(\rho)$, and if this is positive, then the cogood position in $\pi\ast\rho$ corresponds to the cogood position in $\rho$.
\item
Suppose $\phi(\pi)>\ep(\rho)$.  Then $\phi(\pi\ast\rho)=\phi(\pi)+\phi(\rho)-\ep(\rho)$, and the cogood position in $\pi\ast\rho$ coincides with the cogood position in $\pi$.
\end{itemize}
\end{lemma}

\sect{Some $\fsl_3$-crystals}\label{sl3sec}

In this section, we give a construction of rank $2$ crystals, and use Proposition \ref{aaas} to prove that these crystals are regular $\fsl_3$-crystals.  The crystals defined in this way will arise as rank $2$ subcrystals of our $\hsl_n$-crystals, but it will be helpful to describe them in a more abstract way here.

Throughout this section, $i$ and $j$ should be regarded as abstract symbols; in Section \ref{proofsec} they will take values in $\znz$.

\subsection{Biorders}\label{biordsec}

\begin{defn}
Define a \emph{biorder} to be a finite set $S$, equipped with two total orders $>_i$ and $>_j$ and a function $\vh:S\to\{i,j\}$, satisfying the following condition: if $s,t\in S$ are such that $s>_it>_js$, then $\vh(s)=j$ and $\vh(t)=i$.
\end{defn}

Given a biorder $S$ and given $s,t\in S$, we define $s\gg t$ if $s>_it$ and $s>_jt$.  We also define $s\g t$ if either $s>_it$ or $s>_jt$, and we define $\doms$ to be the transitive closure of $\g$.  Then $\gg$ is a partial order on $S$, while $\doms$ is a preorder.

\begin{eg}
Let $S=\{p,q,r,s,t\}$, with
\[p>_iq>_ir>_is>_it,\qquad r>_jp>_jq>_jt>_js.\]
Then we must have
\[\vh(p)=\vh(q)=\vh(s)=j,\qquad \vh(r)=\vh(t)=i.\]
The partial order $\gg$ is given by the Hasse diagram
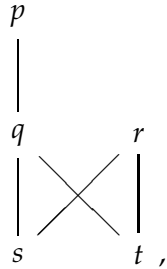
\[\raisebox{45pt}{$\begin{diagram}
\node{p}\arrow{s,-}{}\\
\node{q}\arrow{s,-}{}\arrow{se,-}{}\node{r}\arrow{sw,-}{}\arrow{s,-}{}\\
\node{s}\node{t}
\end{diagram}$},\]
while the preorder $\doms$ is defined by $x\doms y$ if and only if $x\in\{p,q,r\}$ or $y\in\{s,t\}$.
\end{eg}

  If $S$ is a biorder, we define a \emph{configuration} for $S$ to be a function $a:S\to \{0,1,2\}$; we write $\cc(S)$ for the set of configurations of $S$.  Given $a\in\cc(S)$, we shall find it helpful to abuse notation and incorporate the function $\vh$ into $a$; for example, we may write $a(s) = 1_i$ to mean that $a(s)=1$ and $\vh(s)=i$.

We shall define functions $\ee i,\ff i,\ee j,\ff j$ which will make $\cc(S)$ into an $\{i,j\}$-crystal.  Take $a\in\cc(S)$, and for $s\in S$ define
\[\pi_i(a,s)=\begin{cases}
+ & (\text{if $a(s)=0_i$ or $1_j$})\\
- & (\text{if $a(s)=1_i$ or $2_j$})\\
0 & (\text{otherwise}).
\end{cases}\]
Now define the \emph{$i$-signature} of $a$ to be the $\pm$-sequence $\pi_i(a,s_1)\dots\pi_i(a,s_m)$, where $s_1,\dots,s_m$ are the elements of $S$ arranged so that $s_1>_i\dots>_is_m$.

If there is no good position in the $i$-signature of $a$, then define $\ee ia=0$; otherwise, let $k$ be the good position, and define $\ee ia$ to be the configuration with
\[(\ee ia)(s) = a(s)-\delta_{ss_k}.\]
We shall say that $s_k$ is the $i$-good element of $S$ for $a$.

Similarly, if there is no cogood position in the $i$-signature, then we define $\ff ia=0$; otherwise, we let $l$ be the cogood position, and define $\ff ia$ to be the configuration with
\[(\ff ia)(s) = a(s)+\delta_{ss_l};\]
we say that $s_l$ is the $i$-cogood element of $S$ for $a$.

By interchanging the symbols $i$ and $j$ throughout the above definition, we obtain the definition of the $j$-signature of a configuration and functions $\ee j$ and $\ff j$.  The next lemma shows that these functions make $\cc(S)$ into an $\{i,j\}$-crystal.

\begin{lemma}\label{cccrystal}
Suppose $S$ is a biorder, and $h\in\{i,j\}$.
\begin{enumerate}
\item
If $a,b\in\cc(S)$, then $\ee ha=b$ if and only if $\ff hb=a$.
\item
If $a\in\cc(S)$ and we define
\[\ep_h(a)=\max\big\{\ep\ \big|\ \ee h^\ep a\neq0\big\},\qquad \phi_h(g)=\max\big\{\phi\ \big|\ \ff h^\ep a\neq0\big\},\]
then $\ep_h(a)=\ep(\pi)$ and $\phi_h(a)=\phi(\pi)$, where $\pi$ is the $h$-signature of $a$.
\end{enumerate}
\end{lemma}

\begin{pfenum}
\item
We shall suppose that $b=\ee ha$ and prove that $a=\ff hb$; the other direction is similar.  We have
\[b(s) = a(s)-\delta_{ss'}\]
for each $s\in S$, where $s'$ is the $h$-good element of $S$ for $a$.  By Lemma \ref{pmutil}(\ref{pmu2}), we must have $\pi_h(a,s')=-$, and therefore (from the definition of $\pi_h$) we have $\pi_h(b,s')=+$; clearly $\pi_h(b,s)=\pi_h(a,s)$ for all $s\neq s'$.  Now by Lemma \ref{goodcogood} $s'$ is the $h$-cogood element of $S$ for $b$, and therefore $\ff hb=a$.
\item
We prove that $\ep_h(a)=\ep(\pi)$ by induction on $\ep(\pi)$.  If $\ep(\pi)=0$, then by definition $\ee ha=0$, so that $\ep_h(a)=0$.  If $\ep(\pi)>0$, then $c$ has a good position, and so $\ee ha\neq0$.  We write $b=\ee ha$.  From above, the $h$-signature $\rho$ of $b$ is obtained from $\pi$ by replacing $-$ with $+$ in the good position.  Hence by Lemma \ref{goodcogood} we have $\ep(\rho)=\ep(\pi)-1$; by induction $\ep_h(b)=\ep(\rho)$, and the result follows.

Proving that $\phi_h(a)=\phi(\pi)$ is very similar.
\end{pfenum}

\begin{eg}
Let $S=\{p,q,r\}$, with $r>_ip\gg q>_jr$ (so that necessarily $\vh(p)=\vh(q)=i$ and $\vh(r)=j$).  Then the crystal graph of $\cc(S)$ is given in Figure 2 (where we write a configuration $a$ in the form $(a(p),a(q),a(r))$).
\begin{figure}[p]
\setlength\unitlength{1.5pt}
\[\begin{picture}(180,200)
\multiput(120,200)(0,-60)3{\put(0,0){\jaa}\put(60,-40){\iaa}}
\multiput(0,0)(0,80)2{\put(60,30){\ja}\put(60,90){\ja\ia}\put(0,60){\ja}\put(120,60){\ia}\put(120,120){\ia}}
\put(60,110){\vt{2,1,0}}
\put(60,90){\vt{1,0,2}}
\put(0,60){\vt{1,1,2}}
\put(0,140){\vt{1,1,0}}
\put(60,30){\vt{2,1,2}}
\put(60,170){\vt{1,0,0}}
\put(120,0){\vt{2,2,2}}
\put(120,60){\vt{2,0,2}}
\put(120,80){\vt{2,2,0}}
\put(120,120){\vt{0,0,2}}
\put(120,140){\vt{2,0,0}}
\put(120,200){\vt{0,0,0}}
\put(180,40){\vt{2,2,1}}
\put(180,100){\vt{2,0,1}}
\put(180,160){\vt{0,0,1}}
\end{picture}\]

\[\begin{picture}(120,140)
\put(0,60){\ja\vt{0,1,2}}
\put(0,80){\ja\vt{1,2,0}}
\put(0,140){\ja\vt{0,1,0}}
\put(60,30){\ia\vt{0,2,2}}
\put(60,90){\ia\vt{0,1,1}}
\put(60,110){\ia\jaa\vt{0,2,0}}
\put(120,70){\iaa\vt{0,2,1}}
\put(0,0){\vt{1,2,2}}
\put(60,50){\vt{1,2,1}}
\end{picture}
\qquad
\begin{picture}(60,60)
\put(0,30){\ja}
\put(60,60){\ia}
\put(0,30){\vt{1,1,1}}
\put(60,60){\vt{1,0,1}}
\put(60,0){\vt{2,1,1}}
\end{picture}\]

\vspace{\topsep}
\begin{center}Figure 2\end{center}
\end{figure}
\end{eg}

\subsection{Good configurations}

We see from Figure 2 that not every component of $\cc(S)$ is a regular $\fsl_3$-crystal: the second component pictured fails to satisfy axiom A3, and does not have a unique source.  So we need to restrict attention to certain components of $\cc(S)$.

\begin{defn}
Suppose $S$ is a biorder, and $a\in\cc(S)$.  We say that $a$ is \emph{good} for $S$ if it satisfies the following conditions.
\renewcommand{\theenumi}{G\arabic{enumi}}
\begin{enumerate}
\item
There do not exist $s,t\in S$ such that $s>_i t>_j s$ and $a(s)=a(t)=1$.
\item
There do not exist $s,t\in S$ such that $s\g t\doms s$, $\vh(s)=\vh(t)$ and $a(s)<a(t)$.
\item
There do not exist $q,r,s,t\in S$ such that:
\begin{itemize}
\item
$q\g r\gg s\g t$, 
\item
$r\gg t\g q\gg s$,
\item
$a(q)=a(s)=2$, $a(r)=a(t)=0$.
\end{itemize}
\end{enumerate}
\renewcommand{\theenumi}{\arabic{enumi}}
\end{defn}

We write $\gc(S)$ for the set of good configurations for $S$.  We want to show that $\gc(S)$ is a union of connected components of the crystal $\cc(S)$.  That is, we prove the following.

\begin{propn}\label{goodef}
Suppose $S$ is a biorder, $h\in\{i,j\}$ and $a,b$ are configurations for $S$ with $\ff ha=b$.  Then $a$ is good if and only if $b$ is good.
\end{propn}

First we need a lemma concerning axiom G2; it says that if a configuration contains a counterexample to G2, then it contains a counterexample of a particular form.

\begin{lemma}\label{g2tri}
Suppose $S$ is a biorder, and $a$ is a configuration for $S$ which does not satisfy axiom G2.  Then there exist $s,t,r\in S$ such that
\[s\gg t\g r\g s,\quad \vh(s)=\vh(t),\quad a(s)<a(t).\]
\end{lemma}

\begin{pf}
By hypothesis, we can find $m>0$ and $s,t,r_1,\dots,r_m\in S$ such that
\[s\g t\g r_1\g\cdots\g r_m\g s,\qquad \vh(s)=\vh(t),\qquad a(s)<a(t).\]
All we need to do is show that we can make such a choice with $m=1$; the conditions $s\g t$ and $\vh(s)=\vh(t)$ guarantee that $s\gg t$.

So suppose we have chosen $s,t,r_1,\dots,r_m$ as above with $m$ as small as possible, and suppose for a contradiction that $m>1$.  Note first that $r_2\gg t$; for if not, then $s\g t\g r_2\g r_3\g\cdots\g r_m\g s$, contradicting the minimality of $m$.  So we have $r_1\g r_2\gg t\g r_1$, and this implies
\[t\g r_1\g t\qquad\text{and}\qquad r_1\g r_2\g r_1.\]
By the definition of a biorder, we then have $\vh(t)\neq\vh(r_1)\neq\vh(r_2)$, so that $\vh(r_2)=\vh(t)$.

In a similar way, we find that $\vh(r_m)\neq\vh(s)$.  Since $\vh(r_2)=\vh(t)=\vh(s)$, we see in particular that $m>2$.  This then implies that $s\gg r_2$; for if $r_2\g s$, then we have $s\g t\g r_1\g r_2\g s$, contradicting the minimality of $m$.

Since $a(s)<a(t)$, we must have either $a(s)<a(r_2)$ or $a(r_2)<a(t)$.  In the first case, we get
\[s\gg r_2\g r_3\g\cdots\g r_m\g s,\quad \vh(s)=\vh(r_2),\quad a(s)<a(r_2),\]
while in the second case we get
\[r_2\gg t\g r_1\g r_2,\quad \vh(r_2)=\vh(t),\quad a(r_2)<a(t);\]
either way, we have a contradiction to the minimality of $m$.
\end{pf}

\begin{pfofnb}{Proposition \ref{goodef}}
We assume that $h=i$; the case where $h=j$ then follows from the symmetry of the definitions.  We assume that $a$ is not a good configuration, and prove that $b$ is not good; the reverse direction follows in a very similar way.

Let $s'$ be the $i$-cogood element for $a$.  Then $b$ is given by
\[b(s) = a(s)+\delta_{ss'}\]
for $s\in S$, and $a(s')$ equals either $0_i$ or $1_j$.

The assumption that $a\notin\gc(S)$ means that $a$ violates one of the axioms G1--3.  We consider each of the possibilities in turn.
\begin{description}
\item[$a$ does not satisfy G1]
In this case there are $s,t\in S$ such that $s>_i t>_j s$ and $a(s)=a(t)=1$.  From the definition of a biorder, we must have $\vh(s)=j$ and $\vh(t)=i$; hence we cannot have $s'=t$.  If $s'\neq s$, then obviously $b$ does not satisfy G1, so we suppose that $s'=s$.

We have $\pi_i(b,s)=\pi_i(b,t)=-$.  Since $s>_i t$, the $-$ corresponding to $t$ occurs after the $-$ corresponding to $s$ in the $i$-signature of $b$.  But the $-$ corresponding to $s$ is in the good position, so by Lemma \ref{pmutil}(\ref{pmu5}) there must be a $+$ between these two positions.  That is, there is $r\in S$ such that $s>_ir>_it$ and $b(r)$ equals either $0_i$ or $1_j$.

If $b(r)=0_i$, then we have
\[r\g t\g s\g r,\quad \vh(r)=\vh(t),\quad b(r)<b(t),\]
so $b$ violates G2.  On the other hand, if $b(r)=1_j$, then, we have $s\gg r$ by the definition of a biorder, whence
\[r>_i t>_jr,\quad b(r)=b(t)=1,\]
so $b$ violates G1.  Either way, we find $b\notin\gc(S)$.
\item[$a$ does not satisfy G2]
In this case, we apply Lemma \ref{g2tri}, and we find that there are $r,s,t\in S$ such that $s\gg t\g r\g s$, $\vh(s)=\vh(t)$ and $a(s)<a(t)$.  If the same $r,s,t$ do not yield a violation of G2 in $b$, we must be in one of the following two situations:
\begin{enumerate}
\item\label{g21}
$s'=s,\quad a(s)=0_i,\quad a(t)=1_i,\quad \vh(r)=j,\quad t>_jr>_is$;
\item\label{g22}
$s'=s,\quad a(s)=1_j,\quad a(t)=2_j,\quad \vh(r)=i,\quad t>_ir>_js$.
\end{enumerate}
As above, the fact $\pi_i(b,s)=\pi_i(b,t)=-$ while $s>_it$ means that there must be some $q\in S$ such that $s>_iq>_it$ and $b(q)$ equals either $0_i$ or $1_j$.

\begin{itemize}
\item
Suppose we are in case \ref{g21}.  If $b(q)=0_i$ then $b$ violates G2, since
\[q\g t\doms q,\quad \vh(q)=\vh(t),\quad b(q)<b(t).\]
So suppose instead that $b(q)=1_j$.  Then we have $r>_iq$, and since $\vh(q)=j$, this means that $r>_jq$ too.  Hence $t>_jq$, and so we have
\[q>_it>_jq,\quad b(q)=b(t)=1,\]
so $b$ violates G1.
\item
Next suppose we are in case \ref{g22}.  In this case, if we have $b(q)=1_j$, then $b$ violates G2 (via the pair $(q,t)$), so we suppose instead that $b(q)=0_i$.  Now we have $q>_ir$, so
\[q\g r\doms q,\quad \vh(q)=\vh(r),\quad b(q)=0,\]
so that if $b(r)>0$ then $b$ violates G2.  So let us suppose that $b(r)=0$.  Then we find that
\[s\g q\gg t\g r,\qquad q\gg r\g s\gg t,\qquad b(s)=b(t)=2,\qquad b(q)=b(r)=0,\]
and $b$ violates G3.
\end{itemize}
\item[$a$ does not satisfy G3]
Suppose $q,r,s,t$ are as in G3.  Note that the axioms for a biorder imply that $\vh(q)\neq\vh(r)\neq\vh(s)\neq\vh(t)$.  Obviously if $s'\notin\{q,r,s,t\}$, then $b$ violates G3, so suppose otherwise.  $s'$ cannot equal either $q$ or $s$.  If $s'=t$, then $b$ violates G2, because
\[r\g t\doms r,\quad \vh(r)=\vh(t),\quad b(r)<b(t);\]
so suppose $s'=r$.  This means that $\vh(r)=i$ and hence $\vh(s)=j$, so that $\pi_i(b,r)=\pi_i(b,s)=-$.  Arguing as in previous cases, we find that there must be some $p\in S$ such that $r>_ip>_is$ and $b(p)$ equals either $0_i$ or $1_j$.

If $b(p)=1_j$, then $b$ violates G2, since we have
\[p>_i s\doms p,\quad \vh(p)=\vh(s),\quad b(p)<b(s);\]
so suppose instead that $b(p)=0_i$.  Since $p>_i s$ and $\vh(p)=i$, the definition of a biorder implies that $p\gg s$, so we have
\[q\g p\gg s\g t,\qquad p\gg t\g q\gg s,\qquad b(q)=b(s)=2,\qquad b(p)=b(t)=0,\]
and $b$ violates G3.\hfill\qedsymbol
\end{description}
\indent\end{pfofnb}

Given Proposition \ref{goodef}, it makes sense to refer to a connected component of $\cc(S)$ as \emph{good} if every vertex is labelled by a good configuration, and to view $\gc(S)$ as a crystal.  Our aim is to show that every component of $\gc(S)$ is a regular $\fsl_3$-crystal.

\subsection{Axiom A2 holds in $\gc(S)$}

First we note the following simple lemma, which follows from the definitions.

\begin{lemma}\label{esig}
Suppose $a,b\in\cc(S)$ with $b=\ee ia$.  Then the $j$-signature of $b$ is obtained from the $j$-signature of $a$ either by replacing a $0$ with a $-$ or by replacing a $+$ with a $0$.
\end{lemma}

This lemma, together with Lemmata \ref{downbyone} and \ref{cccrystal}, implies that the first part of axiom A2 holds in $\cc(S)$ (and in particular, in $\gc(S)$).  In fact, we can make a more precise statement.

\begin{lemma}\label{epchange}
Suppose $a,b\in\cc(S)$ with $b=\ee ia$.  If $\phi_j(a)=0$, then we have
\[\ep_j(b)=\ep_j(a)+1,\quad \phi_j(b)=\phi_j(a).\]
If not, then let $t$ be the $i$-good element of $S$ for $a$, and let $s$ be the $j$-cogood element for $a$.
\begin{itemize}
\item
If $t>_js$, then
\[\ep_j(b)=\ep_j(a)+1,\qquad \phi_j(b)=\phi_j(a).\]
\item
If $s\gs_jt$, then
\[\ep_j(b)=\ep_j(a),\qquad \phi_j(b)=\phi_j(a)-1,\]
and if in addition $\ep_j(a)>0$ then the $j$-good element for $b$ is the same as the $j$-good element for $a$.
\end{itemize}
\end{lemma}

\begin{pf}
This follows from Lemmata \ref{downbyone}, \ref{cccrystal} and \ref{esig}.
\end{pf}

Of course, the same result holds with $i$ and $j$ interchanged.  Now to complete the proof that axiom A2 holds in $\gc(S)$, we need to show the following.

\begin{lemma}
There is no $a\in\gc(S)$ such that
\[\ee ia,\ff ia\neq0,\quad \ep_j(\ee ia)=\ep_j(a),\quad \phi_j(\ff ia)=\phi_j(a).\]
\end{lemma}

\begin{pf}
Suppose for a contradiction that we can find such a configuration $a$.  Then obviously the $i$-signature of $a$ has both a good and a cogood position.  By Lemma \ref{epchange} the $j$-signature of $a$ has a cogood position, and similarly the $j$-signature must also have a good position.  We let $q,r,s,t$ be the $i$-cogood, $j$-good, $j$-cogood and $i$-good elements for $a$, respectively.

We have $t>_iq$ and $r>_js$, by Lemma \ref{pmutil}(\ref{pmu4}).  We also have $s\gs_jt$, by Lemma \ref{epchange}, and in a symmetrical way we have $q\gs_jr$.

So we have $t>_iq>_jt$, and hence $\vh(q)=i$, $\vh(t)=j$.  Since $\pi_i(a,q)=+$ and $\pi_i(a,t)=-$, we must therefore have $a(q)=0_i$ and $a(t)=2_j$.

Now consider the value of $\vh(r)$.  If $\vh(r)=i$, then (since $\pi_j(a,r)=-$) we have $a(r)=2_i$.  But then we have $q\gg r$ (since $q>_jr$ and $\vh(r)=i$), with $\vh(q)=\vh(r)$, $r\doms q$ and $a(q)<a(r)$, and this contradicts axiom G2.  On the other hand, if we have $\vh(r)=j$, then $a(r)=1_j$.  But then we have $r>_jt$, $\vh(r)=\vh(t)$, $t\doms r$ and $a(r)<a(t)$, which again contradicts axiom G2.
\end{pf}

Of course, the same result holds with $i$ and $j$ interchanged, and we see that axiom A2 holds in $\gc(S)$.

\subsection{Axiom A3 holds in $\gc(S)$}

\begin{lemma}
Suppose $a\in\gc(S)$ with $\ee ia,\ee ja\neq0$ and $\ep_j(\ee ia)=\ep_j(a)$.  Then $\ee i\ee ja=\ee j\ee ia\neq0$, and $\ep_i(\ee ja)=\ep_i(a)+1$.
\end{lemma}

\begin{pf}
Let $r$ and $t$ be the $j$-good and $i$-good elements for $a$ respectively.  By Lemma \ref{epchange}, $r$ is also the $j$-good element for $\ee ia$.  So in order to show that $\ee i\ee ja=\ee j\ee ia\neq0$, we need to show that $t$ is also the $i$-good element for $\ee ja$; for then we shall have $\ee i\ee ja=\ee j\ee ia=b$, where
\[b(s) = a(s)-\delta_{sr}-\delta_{st}.\]


Let $s$ be the $j$-cogood element for $a$.  (Note that there must be such an element, by Lemma \ref{epchange}; moreover, that lemma implies that $s\gs_jt$.  We also have $r>_js$ by Lemma \ref{pmutil}(\ref{pmu4}).)

The $i$-signature of $\ee ja$ is obtained from the $i$-signature of $a$ by subtracting $1$ from the entry corresponding to $r$.  If $r>_it$, then by Lemma \ref{downbyone} the position corresponding to $t$ will still be good in the $i$-signature of $\ee ja$, which is what we require; so suppose otherwise, i.e.\ $t>_ir$.

Now we have $r>_js\gs_jt>_ir$; hence $r>_jt>_ir$, which implies that $\vh(r)=i$ and $\vh(t)=j$.  Since $\pi_j(a,r)=\pi_i(a,t)=-$, we then have $a(r)=a(t)=2$.

Now consider $s$.  Because $\pi_j(a,s)=+$, we have either $a(s) = 1_i$ or $a(s)=0_j$.  In the latter case we would have $s\g t\doms s$, with $\vh(s)=\vh(t)$ and $a(s)<a(t)$, and this contradicts G2.  So instead we must have $a(s)=1_i$.  Hence $\pi_i(a,s)=-$, with $t>_is$; since $t$ is the $i$-good element for $a$, there must therefore be some $q\in S$ such that $t>_iq>_is$ and $\pi_i(a,q)=+$.  We consider the two possibilities for $a(q)$.

If $a(q)=1_j$, then we have $t\gg q$, so that $s>_j q>_is$; but then we have a contradiction to G1.  On the other hand, if $a(q)=0_i$, then we have $q\g s\doms q$ with $\vh(q)=\vh(s)$ and $a(q)<a(s)$, and we have a contradiction to G2.

It remains to show that $\ep_i(\ee ja)=\ep_i(a)+1$.  If this is not the case, then by Lemma \ref{epchange} (with $i$ and $j$ interchanged) we have $\phi_i(a)>0$, and $p\gs_ir$, where $p$ is the $i$-cogood element of $S$.  So we have $r>_js\gs_jt>_ip\gs_ir$; in particular, $r>_jt>_ir$, so that $\vh(r)=i$ and $\vh(t)=j$.  Hence $a(r)=a(t)=2$.

We have $a(p)=0_i$ or $1_j$.  If $a(p)=0_i$, then $p\g r\doms p$ with $\vh(p)=\vh(r)$ and $a(p)<a(r)$, contradicting axiom G2.  So instead $a(p)=1_j$.  Similarly $a(s)=1_i$.  Now we have $\vh(t)=\vh(p)$ and $t\g p$, which implies $t\gg p$, and hence $s>_jp$.  We also have $r\gg s$, which implies $p>_is$.  So we have $p>_is>_jp$ and $a(p)=a(s)=1$, but this contradicts G1.
\end{pf}

The same result holds with $i$ and $j$ interchanged, so part (a) of axiom A3 holds in $\gc(S)$.  Part (b) is proved in the same way.

\subs{The transitive case}

Now suppose that the biorder $S$ is such that $s\doms t$ for every $s\neq t\in S$; we shall say that $S$ is \emph{transitive} if this is the case.  By Proposition \ref{aaas}, in order to show that every component of $\gc(S)$ is regular, it is enough to show that any such component has only one source.  In fact, we shall prove that under the assumption of transitivity $\gc(S)$ has only one source; this implies in particular that it has only one component.

\begin{propn}\label{transource}
Suppose $S$ is a transitive biorder.  Then $\gc(S)$ has only one source.
\end{propn}

\begin{pf}
It is clear that the configuration $a$ with $a(s)=0$ for all $s\in S$ is good and is a source.  So we must show that there is no source $b$ which is good and has $b(s)>0$ for some $s$.

Suppose $b$ is such a configuration, and suppose first that $b(s)=1$.  Without loss of generality, we assume $\vh(s)=i$, and in fact we assume that $s$ is maximal with respect to the order $>_i$ such that $b(s)=1_i$.

The $i$-signature of $b$ contains a $-$ corresponding to $s$.  Since $b$ is a source, this signature cannot contain a good position, and therefore there must be a $+$ preceding this $-$ in the $i$-signature.  Choose such a $+$, and let $t$ be the corresponding element of $S$.  Then $b(t)$ equals either $0_i$ or $1_j$.  In the first case, we have $t>_is$, $\vh(t)=\vh(s)$, $s\doms t$ (because $S$ is transitive) and $b(t)<b(s)$, and this contradicts axiom G2.  So instead we must have $b(t)=1_j$.  If $s>_jt$, then $b$ violates axiom G1, so we must have $t\gg s$.  Repeating the above argument with $t$ in place of $s$ and with $i$ and $j$ interchanged, we find that there is some $t'\in S$ such that $t'\gg t$ and $b(t')=1_i$.  This implies that $t'\gg s$, but this contradicts the choice of $s$.

So we cannot have $b(s)=1$ for any $s$.  We choose some $s$ such that $b(s)=2$.  Without loss of generality, we assume that $\vh(s)=j$, and we assume that $s$ is maximal with respect to the order $>_i$ such that $b(s)=2_j$.  Arguing as above, there must be some $t\in S$ such that $t>_is$ and there is a $+$ corresponding to $t$ in the $i$-signature of $b$.  The assumption that $b(t)\neq1$ now implies that $b(t)=0_i$.  Now the definition of a biorder gives $t\gg s$.

Since $S$ is transitive, we can find $r_1,\dots,r_m\in S$ such that $s\g r_1\g\cdots\g r_m\g t$.  We make such a choice with $m$ as small as possible.  We note first that $m$ must be greater than $1$.  For if $m=1$, then we have $s\g r_1\g s$, which means that $\vh(r_1)=i$; but we also have $r_1\g t\g r_1$, which implies that $\vh(r_1)=j$.

We observe that $t\gg r_{m-1}$.  Indeed, if not, then $s\g r_1\g\cdots\g r_{m-1}\g t$, which contradicts the minimality of $m$.  Now the fact that $r_{m-1}\g r_m$ implies that $t\g r_m$.  So $t\g r_m\g t$, and we therefore have $\vh(r_m)=j$, and $r_m>_it>_jr_m$.  This gives $r_m>_is$, and now the choice of $s$ implies that $b(r_m)=0$.  But then we have $r_m\gg s\doms r_m$, $\vh(r_m)=\vh(s)$ and $b(r_m)<b(s)$, and this contradicts axiom G2.
\end{pf}

Now Proposition \ref{aaas} implies the following.

\begin{cory}\label{transreg}
Suppose $S$ is a transitive biorder.  Then $\gc(S)$ has exactly one component, which is a regular $\fsl_3$-crystal.
\end{cory}

\subs{Every component of $\gc(S)$ is a regular crystal}

Now we complete the proof of the main result of this section by proving that if $S$ is a (possibly intransitive) biorder, then every component of $\gc(S)$ is a regular $\fsl_3$-crystal.  In order to do this, we use tensor products and exploit Proposition \ref{tenshw}.

\begin{propn}\label{intranstens}
Suppose $S$ is a biorder, and suppose that we can decompose $S$ as $S_1\sqcup S_2$ in such a way that for all $s\in S_1$ and $t\in S_2$ we have $s\gg t$.  Then as crystals we have $\cc(S)\cong\cc(S_1)\otimes\cc(S_2)$.  Furthermore, under this isomorphism, every good component of $\cc(S)$ arises as a component of $D_1\otimes D_2$, where $D_1$ is a good component of $\cc(S_1)$ and $D_2$ is a good component of $\cc(S_2)$.
\end{propn}

\begin{pf}
Given a configuration $a$ for $S$ and given $g\in\{1,2\}$, we define $a_g$ to be the restriction of $a$ to $S_g$; so $a_g$ is a configuration for $S_g$.  This defines a bijection
\begin{align*}
\chi: \cc(S)&\lra \cc(S_1)\times\cc(S_2)\\
a&\longmapsto (a_1,a_2),
\end{align*}
and we claim that this bijection gives an isomorphism of crystals, i.e.\ for $a\in\cc(S)$ and $h\in\{i,j\}$
\begin{align*}
\chi(\ee ha) &= \begin{cases}
(\ee ha_1,a_2) & (\phi_h(a_1)\gs\ep_h(a_2))\\
(a_1,\ee ha_2) & (\phi_h(a_1)<\ep_h(a_2)),
\end{cases}\\
\chi(\ff ha) &= \begin{cases}
(\ff ha_1,a_2) & (\phi_h(a_1)>\ep_h(a_2))\\
(a_1,\ff ha_2) & (\phi_h(a_1)\ls\ep_h(a_2))
\end{cases}
\end{align*}
(where we interpret $(a_1,0)$ and $(0,a_2)$ as $0$).  To see this, we note that the $h$-signature of $a$ consists of the $h$-signature of $a_1$ followed by the $h$-signature of $a_2$, and then apply Lemma \ref{concat}.

For the second part of the proposition, we simply observe that if $a$ is a good configuration for $S$, then $a_1$ and $a_2$ are good configurations for $S_1$ and $S_2$ respectively; this is immediate from the definition of a good configuration.
\end{pf}

Finally we can prove the main result of this section.

\begin{thm}\label{mainbiorder}
Suppose $S$ is a biorder.  Then every component of $\gc(S)$ is a regular $\fsl_3$-crystal.
\end{thm}

\begin{pf}
We proceed by induction on $|S|$.  If $S$ is transitive, then the result follows from Corollary \ref{transreg}, so suppose $S$ is not transitive.  This means that we can write $S=S_1\sqcup S_2$, where $S_1,S_2$ are non-empty and $s\gg t$ for all $s\in S_1$ and $t\in S_2$ (for example, choose $s_0,t_0\in S$ such that $t_0\ndom s_0$, and let $S_1=\left\{s\in S\ \left|\ s\dom s_0\right.\right\}$).  By induction every good component of $\cc(S_1)$ or $\cc(S_2)$ is a regular $\fsl_3$-crystal; by Proposition \ref{intranstens}, every good component of $\cc(S)$ is isomorphic to a connected component of the tensor product of a good component of $\cc(S_1)$ and a good component of $\cc(S_2)$, and so by Proposition \ref{tenshw} is a regular $\fsl_3$-crystal.
\end{pf}

\sect{Proof of the main theorem}\label{proofsec}

Now we can return to our crystal $\regg A$ for a given arm sequence $A$, and prove Theorem \ref{main}.  From now on, we write $I$ for the set $\znz$.

\subs{$\regg A$ is an $\hsl_n$-crystal}

The first thing we need to do in order to prove Theorem \ref{main} is to show that the operators $\ee i,\ff i$ actually map $\regg A$ to $\regg A\sqcup\{0\}$.  That is, we need to prove the following proposition.

\begin{propn}\label{regcry}
Suppose $A$ is an arm sequence, $\la\in\regg A$ and $i\in\znz$.
\begin{enumerate}
\item
If $\ff i\la\neq0$, then $\ff i\la\in\regg A$.
\item
If $\ee i\la\neq0$, then $\ee i\la\in\regg A$.
\end{enumerate}
\end{propn}

\begin{pfnb}
We prove (1); the proof of (2) is similar.  Let $\mu=\ff i\la$, and suppose $\mu$ is obtained from $\la$ by adding the $i$-node $(a,c)$.  We suppose for a contradiction that $\mu$ has an illegal hook.  Since $\la$ has no illegal hooks, an illegal hook of $\mu$ must be either the $(b,c)$-hook for some $b<a$, or the $(a,d)$-hook for some $d<c$.  We assume the latter; the other case follows in a similar way (or indeed by the fact that the definitions have a symmetry corresponding to conjugation of partitions).

So we assume that the $(a,d)$-hook of $\mu$ has length $nt$ and arm length $c-d=A_t$, for some $t$.  This means that
\[\la'_{d} = a+nt-A_t-1.\]
Let $b=\la'_d+1$.  Then the node $(b,d)$ is not a node of $\la$, but the node $(b-1,d)$ is.  And in fact $(b,d)$ must be an addable node of $\la$; for if not, then $d>1$ and the node $(b,d-1)$ does not lie in $\la$.  But this then means that the $(a,d-1)$-hook of $\la$ is illegal (with length $nt$ and arm length $A_t$).

Since $c-d=A_t$, we have $(b,d)\cg(a,c)$.  
So $\la$ has addable $i$-nodes $(b,d)\cg(a,c)$, and $(a,c)$ is the $i$-cogood node.  By Lemma \ref{pmutil}(\ref{pmu6}), this means that there is a removable $i$-node $(f,g)$ of $\la$ with $(b,d)\cg(f,g)\cg(a,c)$.  We now consider three cases.
\begin{itemize}
\item
First suppose $f<a$.  The fact that $(f,g)$ and $(a,c)$ have the same residue implies that $g-c+a-f=nu$ for some positive integer $u$.  Since $(f,g)\cg(a,c)$, we have $g-c>A_u$.  Combining this with the fact that $c-d=A_t$, we get $g-d>A_t+A_u$, whence $g-d\gs A_{t+u}$.  On the other hand the fact that $(b,d)\cg(f,g)$ means that $g-d\ls A_{t+u}$.  So $g-d=A_{t+u}$, which means that the $(f,d)$-hook of $\la$ has length $n(t+u)$ and arm length $A_{t+u}$, a contradiction.

\item
Next suppose $a<f<b$.  Now we have $g-d+b-f=nu$ for some positive $u$, and we claim that $g-d=A_u$, which means that the $(f,d)$-hook of $\la$ is illegal.

Now the ordering $(b,d)\cg(f,g)\cg(a,c)$ gives $g-d\ls A_u$, and $c-g\ls A_{t-u}$; using the fact that $c-d=A_t$, we get $g-d\gs A_t-A_{t-u}\gs A_u$.  So $g-d=A_u$.

\item
Finally, we suppose $f>b$.  Now we have $d-g+f-b=nu$ for some positive $u$, and we will show that $d-1-g=A_u$, which means that the $(b,g)$-hook of $\la$ is illegal.

The ordering of the nodes gives $d-g>A_u$ and $c-g\ls A_{t+u}$.  The latter yields $d-g\ls A_{t+u}-A_u$, which is at most $A_u+1$.  So $d-g=A_u+1$, as required.\hfill\qedsymbol
\end{itemize}
\indent\end{pfnb}

Now we can show that $\regg A$ is an $I$-crystal, by checking axiom C1.

\begin{lemma}\label{c1}
Suppose $A$ is an arm sequence, $\la,\mu\in\regg A$ and $i\in I$.  Then $\ee i\la=\mu$ if and only if $\ff i\mu=\la$.
\end{lemma}

\begin{pf}
This follows from the definitions, together with Lemma \ref{goodcogood}.
\end{pf}

Next we need to check that $\regg A$ is an $\hsl_n$-crystal, with the weight function $\wt$ given in \S\ref{models}.  It is immediate from the definition of the weight function that axiom C3 is satisfied.  To verify axiom C4, suppose $\la\in\regg A$ and $i\in I$, and let $\pi$ be the $\pm$-sequence corresponding to the addable and removable $i$-nodes of $\la$.  Then by Lemma \ref{goodcogood}, we have $\ep_i(\la)=\ep(\pi)$ and $\phi_i(\la)=\phi(\pi)$ (cf.\ the proof of Lemma \ref{cccrystal}).  Hence by Lemma \ref{pmutil}(1), $\phi_i(\la)-\ep_i(\la)$ equals the number of addable $i$-nodes of $\la$ minus the number of removable $i$-nodes.  It is well-known (and easy to prove by induction on $|\la|$) that this equals $\lan h_i,\wt(\la)\ran$.  So axiom C4 is satisfied, and $\regg A$ is an $\hsl_n$-crystal.

\subs{$\regg A$ has a unique source}

We have seen that $\regg A$ is an $\hsl_n$-crystal; to show that it is isomorphic to the highest-weight crystal $B(\La_0)$, we verify the hypotheses of Theorem \ref{rank2}.  The first thing we need to check is that the crystal $\regg A$ has a unique source.

\begin{propn}\label{onesource}
Suppose $A$ is an arm sequence.  If $\la\in\regg A$ and $\la\neq\varnothing$, then $\la$ has at least one good node.
\end{propn}

\begin{pf}
For this proof, we introduce a partial order on $\bbn^2$: we put $(a,c)\co(b,d)$ if and only if there are $\gamma,\delta\gs0$ such that $(a,c)\cgs (b+\gamma,d+\delta)$.  It is easy to check that $\gs$ is indeed a partial order, which restricts to the order $\cgs$ on the set of nodes of a given residue.

Now let $(a,c)$ be a node of $\la$ which is maximal with respect to the order $\co$, and let $i$ be the residue of $(a,c)$; then we claim that $\la$ has a good node of residue $i$.  Certainly $(a,c)$ is removable, since we have $(a+1,c)\cso(a,c)$ and $(a,c+1)\cso(a,c)$.  Therefore if there is no good $i$-node, then there must be some addable $i$-node $(b,d)$ of $\la$ such that $(b,d)\cg(a,c)$.  We shall show that this forces $\la$ to have an illegal hook, which gives a contradiction.

We assume that $b>a$; the other case is very similar.  Since $(b,d)$ and $(a,c)$ have the same residue, we can write
\[b-a+c-d=nt\]
with $t$ a positive integer.  Then the fact that $(b,d)\cg(a,c)$ implies that $c-d\ls A_t$.  The node $(b-1,d)$ lies in $\la$ (because $(b,d)$ is an addable node of $\la$), so by maximality we have $(b-1,d)\nco(a,c)$.  Hence $(b-1,d)\nco(a,c+1)$, and therefore $(a,c+1)\cg(b-1,d)$ (since these two nodes have the same residue).  This means that $c+1-d>A_t$; so $c-d=A_t$.  So the $(a,d)$-hook of $\la$ has arm length $A_t$; it has length $c-d+b-a=nt$, and therefore is an illegal hook.
\end{pf}

So we know that the empty partition $\varnothing$ is the unique source of $\regg A$; since $\wt(\varnothing)=\La_0$, all that remains is to check the final condition of Theorem \ref{rank2}, namely that every rank $2$ subcrystal of $\regg A$ is regular.

\subs{$\fsl_2\oplus\fsl_2$-subcrystals}

Suppose $i,j\in I$ with $i\neq j\pm1$, and consider the subcrystal of $\regg A$ given by just the $i$- and $j$-arrows.  This subcrystal is an $\fsl_2\oplus\fsl_2$-crystal, so by Proposition \ref{sl2}, all we need to check is the following.

\begin{propn}\label{sl2check}
Suppose $i,j$ are distinct elements of $\znz$ with $j\neq i\pm1$.  Then the operators $\ee i,\ff i$ on $\regg A$ commute with $\ee j,\ff j$.
\end{propn}

\begin{pf}
This is easy to see from the definitions: since $j\neq i\pm1$, two nodes of residues $i,j$ cannot be adjacent.  Therefore applying $\ee i$ or $\ff i$, which involves adding or removing an $i$-node, cannot affect the set of addable and removable $j$-nodes of a partition.  Hence $\ee i,\ff i$ will commute with $\ee j,\ff j$.
\end{pf}

\subs{$\fsl_3$-subcrystals}

Now we consider the rather more awkward case of $\fsl_3$-subcrystals.  Throughout this subsection, we fix $i\in\znz$, and we set $j=i+1$.  We consider the subcrystal $\regg A^i$ of $\regg A$ obtained by deleting all arrows other than those labelled $i$ or $j$.  $\regg A^i$ is an $\{i,j\}$-crystal; we must prove that each component of this crystal is a regular $\fsl_3$-crystal.

Given a connected component $\cla$ of $\regg A^i$, our aim is to construct a biorder $S$ and an isomorphism $\psi:\cla\lra \cld$, where $\cld$ is a component of $\gc(S)$; then by Theorem \ref{mainbiorder} we shall know that $\cla$ is regular.

Let $\la$ be a partition in $\cla$, and define two partitions $\cladw$ and $\claup$ as follows: $\cladw$ is defined by repeatedly removing removable nodes of residues $i$ and $j$ until there are no more; and $\claup$ is defined by repeatedly adding addable nodes of residues $i$ and $j$ until there are no more.  Note that in general, $\claup$ and $\cladw$ need not lie in $\cla$.  However, since $\cla$ is connected, the definition of $\cladw$ and $\claup$ does not depend on the choice of $\la$.

\begin{eg}
Suppose $n=3$, $i=0$ and $\la=(3,2,1^2)$.  Then we have $\cladw=(3,1)$ and $\claup=(5,3,1^2)$.
\end{eg}

Now define an \emph{$i$-domino} in $\bbn^2$ to be a pair of adjacent nodes of which one has residue $i$ and the other has residue $j$.  Since $j=i+1$, the node of residue $j$ in an $i$-domino must lie either immediately above or immediately to the right of the node of residue $i$; we say that the domino is \emph{vertical} in the first case, and \emph{horizontal} in the second case.

With $\claup$ and $\cladw$ defined as above, the set $[\claup]\setminus[\cladw]$ is a disjoint union of dominoes.  We now define a biorder $S_{\cla}$: as a set, this is the set of dominoes in $[\claup]\setminus[\cladw]$.  The order $>_i$ is defined by taking the $i$-nodes in the dominoes, and using the order $\cg$ on these; that is, for $s,t\in S_{\cla}$,
\[s>_it\qquad\Longleftrightarrow\qquad(\text{the $i$-node in $s$})\cg(\text{the $i$-node in $t$});\]
the order $>_j$ is defined in exactly the same way, using the nodes of residue $j$.  The function $\vh$ is defined by
\[\vh(s) = \begin{cases}
i & (\text{if $s$ is horizontal})\\
j & (\text{if $s$ is vertical}).
\end{cases}\]
It is straightforward to check that $S_{\cla}$ satisfies the axioms for a biorder.

Now we define a map $\psi:\cla\rightarrow\cc(S_{\cla})$.  Given a partition $\la$ in $\cla$, we define $\psi(\la)$ by
\[\psi(\la)(s) = \big|[\la]\cap s\big|\]
for each $s\in S_{\cla}$.

\begin{eg}
Continuing the last example, $[\claup]\setminus[\cladw]$ consists of three dominoes:
\begin{itemize}
\item
a horizontal domino $p=\{(1,4),(1,5)\}$;
\item
a horizontal domino $q=\{(2,2),(2,3)\}$;
\item
a vertical domino $r=\{(3,1),(4,1)\}$.
\end{itemize}
If $A$ is an arm sequence with $A_1=1$ and $A_2=3$, then we find that
\[(1,5)\cg(2,3)\cg(3,1),\qquad (4,1)\cg(1,4)\cg(2,2),\]
so that $S_{\cla}$ is precisely the biorder given in the example in \S\ref{biordsec}.  The component $\cla$ of $\regg A^i$ is given in Figure 3 (which the reader should compare with the first diagram in Figure 2); in the Young diagrams, we have marked with $\times$ the nodes belonging to $\cladw$, so that the reader can more easily identify the added nodes.
\begin{figure}
\setlength\unitlength{2.3pt}
\renewcommand{\ja}{\put(10,-5){\vector(2,-1){40}}\put(30,-13.5){\small$1$}}
\renewcommand{\ia}{\put(-10,-5){\vector(-2,-1){40}}\put(-30,-13.5){\small$0$}}
\renewcommand{\jaa}{\put(10,-6,667){\vector(3,-2){40}}\put(30,-18){\small$1$}}
\renewcommand{\iaa}{\put(-10,-6,667){\vector(-3,-2){40}}\put(-30,-18){\small$0$}}
\newcommand{\ts}{\raisebox{-1pt}{$\times$}}
\Yboxdim{8pt}
\[\begin{picture}(180,200)
\multiput(120,200)(0,-60)3{\put(0,0){\jaa}\put(60,-40){\iaa}}
\multiput(0,0)(0,80)2{\put(60,30){\ja}\put(60,90){\ja\ia}\put(0,60){\ja}\put(120,60){\ia}\put(120,120){\ia}}
\put(60,110){\makebox[0pt]{\raisebox{-6pt}{\young(\ts\ts\ts\ \ ,\ts\ )}}}
\put(60,90){\makebox[0pt]{\raisebox{-14pt}{\young(\ts\ts\ts\ ,\ts,\ ,\ )}}}
\put(0,60){\makebox[0pt]{\raisebox{-14pt}{\young(\ts\ts\ts\ ,\ts\ ,\ ,\ )}}}
\put(0,140){\makebox[0pt]{\raisebox{-6pt}{\young(\ts\ts\ts\ ,\ts\ )}}}
\put(60,30){\makebox[0pt]{\raisebox{-14pt}{\young(\ts\ts\ts\ \ ,\ts\ ,\ ,\ )}}}
\put(60,170){\makebox[0pt]{\raisebox{-6pt}{\young(\ts\ts\ts\ ,\ts)}}}
\put(120,0){\makebox[0pt]{\raisebox{-14pt}{\young(\ts\ts\ts\ \ ,\ts\ \ ,\ ,\ )}}}
\put(120,60){\makebox[0pt]{\raisebox{-14pt}{\young(\ts\ts\ts\ \ ,\ts,\ ,\ )}}}
\put(120,80){\makebox[0pt]{\raisebox{-6pt}{\young(\ts\ts\ts\ \ ,\ts\ \ )}}}
\put(120,120){\makebox[0pt]{\raisebox{-14pt}{\young(\ts\ts\ts,\ts,\ ,\ )}}}
\put(120,140){\makebox[0pt]{\raisebox{-6pt}{\young(\ts\ts\ts\ \ ,\ts)}}}
\put(120,200){\makebox[0pt]{\raisebox{-6pt}{\young(\ts\ts\ts,\ts)}}}
\put(180,40){\makebox[0pt]{\raisebox{-10pt}{\young(\ts\ts\ts\ \ ,\ts\ \ ,\ )}}}
\put(180,100){\makebox[0pt]{\raisebox{-10pt}{\young(\ts\ts\ts\ \ ,\ts,\ )}}}
\put(180,160){\makebox[0pt]{\raisebox{-10pt}{\young(\ts\ts\ts,\ts,\ )}}}
\end{picture}\]

\vspace{\topsep}
\begin{center}\normalsize Figure 3\end{center}
\end{figure}
\end{eg}

\begin{lemma}\label{psigood}
Suppose $\la\in\cla$.  Then $\psi(\la)$ is a good configuration.
\end{lemma}

\begin{pfnb}
We use proof by contradiction, showing that if $\psi(\la)$ violates one of axioms G1--3, then $\la$ possesses an illegal hook, so does not lie in $\regg A$.  For brevity, let us write $a$ for $\psi(\la)$.
\begin{description}
\item[$a$ does not satisfy G1]
In this case, we have $s,t\in S_{\cla}$ such that $s>_it>_js$ and $a(s)=a(t)=1$.  This implies that $s$ is a vertical domino, while $t$ is a horizontal domino.  Let $(b,c)$ be the $i$-node in $s$, and $(d,e)$ the $i$-node in $t$.  Then we have $(b,c)\cg (d,e)$, but $(d,e+1)\cg (b-1,c)$.  Since $(b,c)$ and $(d,e)$ have the same residue, we have
\[|b-d+e-c|=nu\]
for some positive integer $u$.

Suppose first that $b>d$.  Now $(b,c)\cg (d,e)$ implies that $e-c\ls A_u$, while $(d,e+1)\cg (b-1,c)$ implies that $e+1-c>A_u$, so we have $e-c=A_u$.  Now consider the $(d,c)$-hook of $\la$.  Since $a(t)=1$, we have $(d,e)\in\la\notni(d,e+1)$, so the arm length of the $(d,c)$-hook is $e-c=A_u$.  Since $a(s)=1$, we have $(b,c)\notin\la\ni(b-1,c)$, and we find that the length of the $(d,c)$-hook is $b-d+e-c=nu$.  So $\la$ possesses an illegal hook.

Alternatively, suppose $d>b$.  Now we have $c-e-1=A_u$, and we claim that the $(b,e)$-hook of $\la$ is illegal.  Since $a(s)=1$ we have $(b,c)\notin\la$.  On the other hand, the definition of $\claup$ implies that $(b,c-1)\in\la$; for otherwise $(b,c)$ could not lie in $\claup$.  So the arm length of the $(b,e)$-hook of $\la$ is $c-1-e=A_u$.  To compute the length of the $(b,e)$-hook, we observe that $(d,e)\in\la$ (because $a(t)=1$), but $(d+1,e)\notin\la$ (by the definition of $\cladw$).  So we find that the length of the $(b,e)$-hook is $d-b+c-e=nu$, as required.
\item[$a$ does not satisfy G2]
Applying Lemma \ref{g2tri}, there are $r,s,t\in S_{\cla}$ such that $s\gg t\g r\g s$, $\vh(s)=\vh(t)$ and $a(s)<a(t)$.  We let $(b,c)$ be the $i$-node in $s$ and $(d,e)$-the $i$-node in $t$.  Now we have $(b,c)\cg(d,e)$, and we claim that $(d+1,e+1)\cg(b,c)$.  To show this, we'll assume $\vh(s)=i$ (the other case being very similar).  In this case, if we let $(f,g)$ be the $i$-node in $r$, then $t\g_j r$ gives $(d,e+1)\cg(f-1,g)$, and hence $(d+1,e+1)\cg(f,g)$.  But we also have $(f,g)\cg(b,c)$ since $r\g_i s$, so $(d+1,e+1)\cg(b,c)$.  
As above, we have $|b-d+e-c|=nu$ for some positive $u$.

First suppose $b>d$.  Then we compute $e-c=A_u$.  This enables us to find an illegal hook in $\la$, but there are various cases.  If $s$ and $t$ are horizontal dominoes, then the $(d,c)$-hook is illegal if $a(s)=0$ and $a(t)=1$, while the $(d,c+1)$-hook is illegal if $a(s)<a(t)=2$.  On the other hand if $s$ and $t$ are vertical, then the $(d,c)$-hook is illegal if $a(s)=1$ and $a(t)=2$, while the $(d-1,c)$-hook is illegal if $a(s)=0<a(t)$.

The case where $d>b$ is similar.
\item[$a$ does not satisfy G3]
Suppose $q,r,s,t\in S_{\cla}$ are as in G3.  Then $\vh(q)\neq\vh(r)\neq\vh(s)\neq\vh(t)$, and we consider the case where $\vh(q)=j$; the case where $\vh(q)=i$ is very similar.  Let $(b,c)$ be the $i$-node in $r$, and $(d,e)$ the $i$-node in $s$.  Then we have $(b,c)\cg (d,e)$, but (since $s>_it>_jq>_ir$) $(d+1,e+1)\cg (b,c)$.  Now, using similar arguments to the previous cases, we find that the $(d,c)$-hook of $\la$ is illegal if $b>d$, while the $(b,e)$-hook is illegal if $d>b$.\hfill\qedsymbol
\end{description}
\indent\end{pfnb}

Finally, we note the following.

\begin{lemma}\label{morph}
$\psi$ commutes with $\ee i,\ff i,\ee j,\ff j$.
\end{lemma}

\begin{pf}
We consider $\ee i$ and $\ff i$; the proof for $\ee j$ and $\ff j$ is similar.  Suppose $\la\in\cla$, and let $\ar i(\la)$ denote the set of addable and removable $i$-nodes of $\la$.  The definitions of $\claup$ and $\cladw$ imply that each node in $\ar i(\la)$ is contained in $[\claup]\setminus[\cladw]$, and so is contained in some domino in $[\claup]\setminus[\cladw]$.  So we have a map $\partial:\ar i(\la)\to S_\cla$, given by mapping a node to the domino that contains it; since each domino contains a unique $i$-node, $\partial$ is injective; furthermore, $\partial$ is order-preserving in the sense that if $(a,c)\cg(b,d)$ in $\ar i(\la)$, then $\partial((a,c))>_i\partial((b,d))$.  Now let $a=\psi(\la)$, and for $s\in S_\cla$ define $\pi_i(a,s)$ as in \S\ref{biordsec}.  Then it is straightforward to check that
\[\pi_i(a,s) = \begin{cases}
+ & (\text{if $s=\partial((a,c))$, with $((a,c))$ an addable node of $\la$})\\
- & (\text{if $s=\partial((a,c))$, with $((a,c))$ a removable node of $\la$})\\
0 & (\text{if $s$ does not lie in the image of $\partial$}).
\end{cases}\]
This, together with the fact that $\partial$ is order-preserving, means that the $\pm$-sequence obtained from $\ar i(\la)$ is precisely the $i$-signature of $\psi(\la)$ with the $0$s removed.  Moreover, if $\pi_i(a,s)=+$, then increasing $a(s)$ by $1$ corresponds to adding the $i$-node in $s$ to $\la$, while if $\pi_i(a,s)=-$, then decreasing $a(s)$ by $1$ corresponds to removing the $i$-node in $s$ from $\la$.  In view of Remark 1 in Section \ref{pmsec}, we see that the definitions of $\ee i$ and $\ff i$ on $\la$ and on $a$ are essentially the same.
\end{pf}

Since $\cla$ is a connected component of $\regg A$, Lemma \ref{morph} implies that the image of $\psi$ is a connected component of $\cc(S)$; this component is good, by Lemma \ref{psigood}.  It is easy to see that $\psi$ is injective, and so $\cla$ is isomorphic to a component of $\gc(S)$, and so by Theorem \ref{mainbiorder} $\cla$ is a regular $\fsl_3$-crystal.

So we have verified all the hypotheses of Theorem \ref{rank2} for the crystal $\regg A$, and the proof of Theorem \ref{main} is complete.

\section{Crystals corresponding to different arm sequences are different}\label{diffsec}

For each $n\gs3$, we have defined a family of crystals for $U_q(\hsl_n)$; or rather, a family of models for the same crystal.  However, it is not clear that different arm sequences give different models.  The aim of this section is to prove this statement, by showing the following.

\begin{propn}\label{different}
Suppose $n\gs3$ and $A,A'$ are distinct arm sequences.  Then the sets $\regg A$ and $\regg{A'}$ are distinct.
\end{propn}

We begin by singling out two particular arm sequences.  Call the two sequences $(0,1,2,\dots)$ and $(n-1,2n-2,3n-3,\dots)$ the \emph{extreme} sequences.

\begin{lemma}\label{nonext}
Suppose $A$ is a non-extreme arm sequence.
\begin{enumerate}
\item
For each $t\gs1$, we have $t\ls A_t\ls (n-1)t-1$.
\item
Let $u\gs1$ and let $\la$ be the partition $(A_u+1,1^{nu-A_u-1})$.  Then $\la$ possesses an illegal hook of length $nu$, but does not possess an illegal hook of any other length.
\end{enumerate}
\end{lemma}

\begin{pfenum}
\item
This is simple to check.
\item
The $(1,1)$-hook is illegal of length $nu$.  Now suppose the $(a,c)$-hook is an illegal hook of length $nt$ for some $t<u$.  Then either $a=1$ and $2\ls c\ls A_u+1$, or $c=1$ and $2\ls a\ls nu-A_u$.  In the first case we find $nt=A_t+1$, which contradicts the fact that $A_t\ls (n-1)t-1$, while in the second case we find that $A_t=0$, which also contradicts the first part of this lemma.
\end{pfenum}

\begin{pfof}{Proposition \ref{different}}
Suppose first that $A$ is the extreme sequence $(n-1,2n-2,\dots)$.  Then $A_1=n-1$, while $A'_1<n-1$, so the partition $(n)$ is $A'$-regular but not $A$-regular.  A similar argument applies in the case where $A=(0,1,2,\dots)$, using the partition $(1^n)$.

Now assume that neither $A$ nor $A'$ is extreme, and let $u$ be minimal such that $A_u\neq A'_u$.  Consider the partition $\la=(A_u+1,1^{nu-A_u-1})$.  By Lemma \ref{nonext}(2), $\la$ does not lie in $\regg A$; nor does it have a hook with length $nt$ and arm length $A_t=A'_t$, for any $t<u$.  The only hook of length $nt$ with $t\gs u$ is the $(1,1)$-hook, which has length $nu$ and arm length $A_u\neq A'_u$, and therefore $\la$ does lie in $\regg{A'}$.
\end{pfof}

\begin{rmks}\indent
\begin{enumerate}
\vspace{-\topsep}\item
The proof of Proposition \ref{different} relies on the restriction $t-1\ls A_t\ls (n-1)t$.  In fact, if we broaden the definition of an arm sequence to allow $0\ls A_t\ls nt-1$, then Theorem \ref{main} still holds.  But the crystals $\regg A$ and $\regg{A'}$ will be identical whenever $A_1=A'_1=0$ or $A_1=A'_1=n-1$.
\item
We can extend the crystal operators to give the structure of a crystal on the set of all partitions (not just the $A$-regular ones).  Then different arm sequences give different crystals, even with the broader definition above.  However, except in very special cases, the components of this crystal other than the component $\regg A$ are not regular.  It seems that some modification of the definitions is appropriate to make these components regular; the author hopes to address this in the future.
\end{enumerate}
\end{rmks}

We end the paper by giving an alternative parametrisation of our crystals; this shows in particular that we have uncountably many arm sequences.  

\begin{lemma}\label{cvgs}
Suppose $A$ is an arm sequence.  Then the sequence $\left(\dfrac{A_t}t\right)$ converges to some $y_A\in[1,n-1]$.
\end{lemma}

\begin{pfnb}
We claim that for any $t,u$ we have
\[\left|\frac{A_t}t-\frac{A_u}u\right|<\frac1{\min\{t,u\}},\]
so that the sequence is Cauchy.  Applying the axioms for an arm sequence repeatedly, we get $uA_t\ls A_{tu}<uA_t+u$ and similarly $tA_u\ls A_{tu}<tA_u+t$.  Hence
\[\raisebox{6pt}{${\displaystyle\frac{A_t}t-\frac{A_u}u< \frac{A_{tu}}{tu}-\frac{A_{tu}-t}{tu}=\frac1u,\qquad\qquad
\frac{A_u}u-\frac{A_t}t< \frac{A_{tu}}{tu}-\frac{A_{tu}-u}{tu}=\frac1t}$}.\tag*{\qedsymbol}\]
\end{pfnb}

We shall abuse terminology by saying that the arm sequence $A$ \emph{tends to} $y$ if $y=\lim_{t\to\infty}\frac{A_t}t$.

\begin{lemma}\label{uniquearm}
Suppose $y\in[1,n-1]$.
\begin{enumerate}
\item
If $y$ is irrational, then the sequence $A^y$ given by
\[A^y_t = \lfloor yt\rfloor\]
is an arm sequence.  Furthermore, $A^y$ is the unique arm sequence that tends to $y$.
\item
If $y$ is rational, then the two sequences $A^{y,+}$ and $A^{y,-}$ given by
\[A^{y,+}_t = \lfloor yt\rfloor,\qquad A^{y,-}_t = \lceil yt-1\rceil\]
are arm sequences.  They are the only arm sequences that tend to $y$.
\end{enumerate}
\end{lemma}

\begin{pf}
It is straightforward to verify that the given sequences are arm sequences tending to $y$.  Now suppose $A$ is an arm sequence tending to $y$.  First we claim that $A_t\ls yt$ for each $t$.  If not, then for some $t$ we have $A_t=yt+\delta$ for $\delta>0$.  For any $N>0$ we get $A_{Nt}\gs NA_t=yNt+N\delta$, and hence
\[\frac{A_{Nt}}{Nt}\gs y+\frac{\delta}t,\]
so $A$ does not tend to $y$; contradiction.

In a similar way, using the inequality $A_{Nt}\ls NA_t+N-1$, we get $A_t\gs yt-1$ for all $t$.  Since each $A_t$ is an integer, this specifies each $A_t$ uniquely, except when $yt$ is an integer.  So the proof of (1) is complete.  To complete the proof of (2), we must show that either $A_t=yt$ whenever $yt\in\bbn$, or $A_t=yt-1$ whenever $yt\in\bbn$; in other words, it is not possible to find $t$ and $u$ such that $A_t=yt-1$ while $A_u=yu$.  But if we do have such $t,u$, then the inequality at the end of the proof of Lemma \ref{cvgs} is violated.
\end{pf}

These results imply that we may parametrise arm sequences by real numbers in the interval $[1,n-1]$, but with two sequences for each rational number in this interval.  This implies that there are uncountably many arm sequences (so we do have uncountably many models of the crystal $B(\La_0)$), and that there is a natural total order on arm sequences.  In the notation of Lemma \ref{uniquearm}, the extreme sequences $(0,1,2,\dots)$ and $(n-1,2n-2,\dots)$ are the sequences $A^{1,-}$ and $A^{n-1,+}$ respectively, and the sequence $(1,2,3,\dots)$ defining Berg's ladder crystal is $A^{1,+}$.

\end{document}